\documentclass[10pt]{article}

\usepackage[
hyperindex  = {true},
colorlinks  = {true},
linkcolor   = {blue},
citecolor   = {blue}
]{hyperref}
\usepackage[sf,bf,medium]{titlesec}
\usepackage{multirow}
\usepackage{color, colortbl}
\usepackage[font=footnotesize]{subcaption}
\usepackage[font=footnotesize]{subcaption}
\usepackage[top=1in,bottom=1.25in,left=1.25in,right=1.25in]{geometry}
\usepackage{ragged2e}
\usepackage{tabto}
\usepackage{latexsym,graphicx, amssymb, amsfonts, setspace, amsmath,amsthm}
\usepackage{color}
\usepackage[capitalise,noabbrev]{cleveref}
\usepackage{marginnote}
\usepackage{lmodern}
\usepackage{multicol}
\usepackage{afterpage}
\usepackage{pgfplots}
\pgfplotsset{compat=1.8}
\usepgfplotslibrary{statistics}
\usepackage{grffile}
\usepackage{pgfplotstable}
\usepackage{slashbox}
\usepackage{mathtools}
\usepackage[utf8]{inputenc}
\usepackage[english]{babel}
\usepackage{float}
\usepackage{esvect}
\usepackage{algorithm,algorithmic}
\usepackage{enumitem}
\usepackage{hhline}
\usepackage{varwidth}
\usepackage{tikz}
\usetikzlibrary{backgrounds,patterns,calc}
\usepackage{xparse}
\usepackage{tabularx}
\usepackage{placeins}

\usepackage{bm}
\usepackage{dsfont}
\usepackage[sort,compress]{cite}
\numberwithin{equation}{section}

\pgfplotsset{compat=newest}

\makeatletter
\newcommand{\doublewidetilde}[1]{{%
  \mathpalette\double@widetilde{#1}%
}}
\newcommand{\double@widetilde}[2]{%
  \sbox\z@{$\m@th#1\widetilde{#2}$}%
  \ht\z@=.9\ht\z@
  \widetilde{\box\z@}%
}
\makeatother

\newtheorem{remark}{\sffamily Remark}[section]
\newtheorem{lemma}{\sffamily Lemma}[section]

\newcommand{\pa}{\partial}
\newcommand{\xb}{\boldsymbol{x}}
\newcommand{\ub}{\boldsymbol{u}}

\newcommand{\db}{\boldsymbol{d}}

\newcommand{\bbR}{{\mathbb R}}
\newcommand{\tran}{^{\top\kern-\scriptspace}}%tranpose symbol
\newcommand{\bnu}{\boldsymbol{\nu}}

\newcommand{\be}{\begin{equation}}
\newcommand{\ee}{\end{equation}}
\newcommand{\ba}{\begin{aligned}}
\newcommand{\ea}{\end{aligned}}
\newcommand{\bea}{\begin{eqnarray}}
\newcommand{\eea}{\end{eqnarray}}

\newcommand\nc{\newcommand}

\newcommand{\bE}{\mathbf{E}}
\newcommand{\bH}{\mathbf{H}}

\newcommand{\uscat}{u^\textrm{scat}}
\newcommand{\uinc}{u^\textrm{inc}}

\newcommand{\cA}{\mathcal{A}}
\newcommand{\cS}{\mathcal{S}}
\newcommand{\cD}{\mathcal{D}}
\newcommand{\cF}{\mathcal{F}}

\nc\ex{E_x}
\nc\ey{E_y}
\nc\ez{E_z}

\nc\hx{H_x}
\nc\hy{H_y}
\nc\hz{H_z}

\nc\px[1]{\frac{\partial #1}{\partial x}}
\nc\py[1]{\frac{\partial #1}{\partial y}}

\nc\br{{\boldsymbol{r}}}
\nc\curl{\nabla\times}
\nc\dive{\nabla\cdot}
\nc\pt{\frac{\partial}{\partial t}}
\nc\pet{\frac{\partial \bE}{\partial t}}
\nc\pht{\frac{\partial \bH}{\partial t}}

\newcommand{\Lim}[1]{\raisebox{0.5ex}{\scalebox{0.8}{$\displaystyle \lim_{#1}\;$}}}

\newcommand{\kmax}{k_{\textrm{max}}}
\newcommand{\bx}{\boldsymbol{x}}
\newcommand{\by}{\boldsymbol{y}}

\DeclarePairedDelimiter{\floor}{\lfloor}{\rfloor}
\DeclareMathOperator*{\argmin}{arg\,min}

\definecolor{Gray}{gray}{0.9}
\newcolumntype{g}{>{\columncolor{Gray}}c}

\pgfmathdeclarefunction{gauss}{3}{%
  \pgfmathparse{1/(#3*sqrt(2*pi))*exp(-((#1-#2)^2)/(2*#3^2))}%
}

\begin{document}

\begin{titlepage}

  \raggedleft
  {\texttt{Technical Report\\
    \today}}
  
  \hrulefill

  \vspace{4\baselineskip}

  \raggedright
  {\LARGE \sffamily\bfseries Multifrequency inverse obstacle scattering 
with unknown impedance boundary conditions using recursive linearization}
  
%  \vspace{\baselineskip}
%  {\Large \sffamily Surface integral equations}

  \vspace{3\baselineskip}
  % author 1
 \vspace{\baselineskip}
 
  \normalsize Carlos Borges\\%\footnote{Enter grant info here.}\\
  \small \emph{Department of Mathematics, University of Central Florida\\
    Orlando, FL 32816}\\
  \texttt{carlos.borges@ucf.edu}
  
  \vspace{\baselineskip}
 
  \normalsize Manas Rachh\\%\footnote{Corresponding author.}\\
  \small \emph{Center for Computational Mathematics, Flatiron Institute\\
    New York, NY 10010}\\
  \texttt{mrachh@flatironinstitute.org}
  \normalsize

\end{titlepage}

\begin{abstract}
In this paper, we consider the reconstruction of the shape and the impedance function of an obstacle from 
measurements of the scattered field at a collection of receivers outside the object. The data is assumed to be generated by plane waves 
impinging on the unknown obstacle from multiple directions and at multiple frequencies. 
This inverse problem can be reformulated as an optimization problem: that of finding band-limited shape and impedance functions
which minimize the $L^2$ distance between the computed value of the scattered field at the receivers and the given measurement data. 
The optimization problem is highly non-linear, non-convex, and ill-posed. Moreover, the objective function is computationally expensive to evaluate (since a large number of Helmholtz boundary value problems need to be solved at every iteration in the optimization loop). The recursive linearization approach (RLA) proposed by Chen has been successful in addressing these issues in the context of recovering the sound speed of an inhomogeneous object or the shape of a sound-soft obstacle. We present an extension of the RLA for the recovery of both the shape and impedance functions of the object. 
The RLA is, in essence, a continuation method in frequency where a sequence of single frequency inverse problems is solved. At each higher frequency, one attempts to recover incrementally higher resolution features using a step assumed to be small enough to ensure that the initial guess
obtained at the preceding frequency lies in the basin of attraction for 
Newton's method at the new frequency. 
We demonstrate the effectiveness of this approach with several numerical examples. Surprisingly, we find that one can recover the shape with high accuracy even when the  measurements are generated by sound-hard or sound-soft objects, eliminating the need to know the precise boundary conditions appropriate for modeling the object under consideration. While the method is effective in obtaining high quality reconstructions for many complicated geometries and impedance functions, a number of interesting open questions remain regarding the convergence behavior of the approach.
We present numerical experiments that suggest underlying mechanisms of
success and failure, pointing out areas where improvements could help lead to robust and automatic tools for the solution of inverse obstacle scattering problems.

\end{abstract}

\small

\tableofcontents

\normalsize
\newpage

\section{Introduction}\label{s:intro}
Scattering problems arise naturally in multiple applications including
medical imaging\cite{kuchment2014radon}, non-destructive testing\cite{collins1995nondestructive,engl2012inverse}, remote sensing \cite{Ustinov2014}, ocean acoustics \cite{chavent2012inverse,0266-5611-10-5-003}, sonar and radar\cite{cheney2009fundamentals}. In many of these applications, a typical inverse problem is finding the shape $\partial \Omega$ of an obstacle  $\Omega$ given far-field or distant measurements of the scattered field. 
In order to formulate a well-posed forward problem,
the obstacle is often assumed to satisfy impedance or Fourier-Robin type
boundary conditions. These can model
complicated wave/surface interactions such as thin coatings \cite{aslanyurek2011generalized}, corrugated or rough surfaces \cite{hu2015variational,lechleiter2009generalized}, or highly absorbing media \cite{haddar2005generalized,haddar2008generalized,nguyen2015generalized}. 
In the time harmonic setting, the forward problem for the scattered field $\uscat$ is given by
\begin{equation} \label{eq:imp_problem}
\begin{cases}
\Delta \uscat+k^2 \uscat = 0, \quad \text{in} ~\Omega, \\
\frac{\partial \uscat}{\partial \nu} + ik \lambda \uscat = -\left(\frac{\partial \uinc}{\partial \nu} + ik \lambda \uinc \right) \, \quad \text{on} ~\pa \Omega \, , \\%= -\left(\bnu \cdot \db + ik \lambda \right) e^{ik \bx \cdot \db} \, \quad \text{on} ~\pa \Omega \, , \\
\Lim{r\rightarrow \infty} r^{1/2}\left(\frac{\partial u^\emph{scat}}{\partial r} - iku^\emph{scat}\right) = 0 \,, 
\end{cases} 
\end{equation}
where $k$ is the wave number, $\lambda$ is the impedance function, $\bnu$ is the normal to the boundary $\partial \Omega$, and $\uinc$ is the incident field. In this work, without loss of generality, we will consider $\uinc = e^{ik \bx \cdot \db}$, a plane wave with wavenumber $k$ and incident direction $\db$. %In the inverse problem, one is interested in recovering the shape of the obstacle $\partial \Omega$ and the impedance function $\lambda$ defined on $\partial \Omega$.% instead of the density and sound speed inside the object $\Omega$.

\begin{figure}[h]
\begin{subfigure}[t]{0.42\textwidth}
\center
\includegraphics[width=1\textwidth]{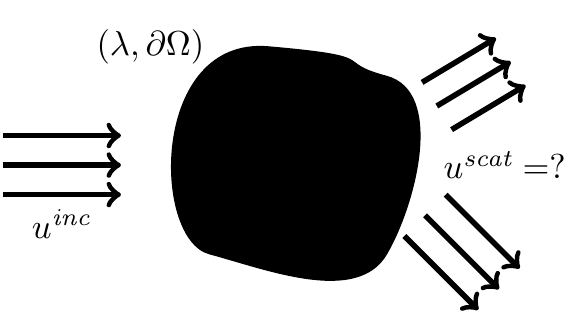}
\caption{Forward scattering impedance problem}\label{fig:fwd_imp_prob}
\end{subfigure}
\hspace{0.2\textwidth}
\begin{subfigure}[t]{0.42\textwidth}
\center
\includegraphics[width=1\textwidth]{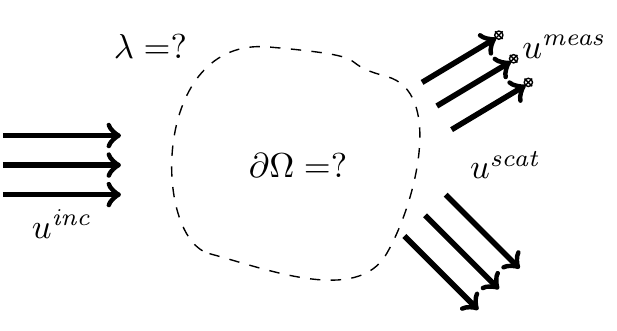}
\caption{Inverse scattering impedance problem}\label{fig:inv_imp_prob}
\end{subfigure}
\caption{Scattering from obstacles. In the {\em forward scattering problem with classical impedance boundary condition}, the shape of the obstacle $\partial \Omega$ and the Impedance function $\lambda$ are given and we want to evaluate the scattered field $u^\emph{scat}$ given the incident field $u^\emph{inc}$, as shown in Figure \ref{fig:fwd_imp_prob}. Finally, in the {\em inverse scattering problem with classical impedance boundary condition}, one wants to recover the shape and impedance function for the obstacle given measurements of the scattered field at receivers place around the obstacle, see Figure \ref{fig:inv_imp_prob}.}\label{fig:problems_intro}
\end{figure}

Given a collection of receivers $\{ \br_{j} \}$, $j=1,2\ldots N_{r}$, we define the forward scattering operator, $\mathcal{F}_{k,\db}:(\partial \Omega,\lambda) \rightarrow \mathbb{C}^{N_{r}}$, as the scattered field evaluated at the receivers, i.e.,
\begin{equation} \label{eq:fwd_imp_op}
\mathcal{F}_{k,\db}(\partial \Omega,\lambda)=\ub_{k,\db}^\emph{meas}
\end{equation}
where the $j$th component of $\ub_{k,\db}^\emph{meas}$ is given by $\uscat_{k,\db}(\br_{j})$. 
The inverse problem corresponding to the forward scattering problem \eqref{eq:fwd_imp_op} is to obtain reconstructions of the shape $\Gamma$ and impedance function $\lambda$, given measurements of the scattered field at the receivers from one or more incident waves at possibly multiple frequencies, see Figure~\ref{fig:problems_intro}.  In particular, given a collection of wavenumbers $k_{j}$, $j=1,2\ldots M$, with $k_{j}>k_{j-1}$, incident directions $\db_{\ell}$, $\ell=1,2,N_{d}$, and the corresponding far field data $\ub_{k_{j},\db_{\ell}}^{\emph{meas}}$, 
the inverse scattering problem seeks the optimum shape $\Gamma$, and impedance function $\lambda$ which minimize the following objective function 
\begin{equation}
\label{eq:multi_freq_inv}
[\tilde{\Gamma},\tilde{\lambda}] = \argmin_{\Gamma,\lambda} \sum_{m=1}^{M} \sum_{\ell=1}^{N_{d}} \| \ub_{k_{m},\db_{\ell}}^{\emph{meas}} - \mathcal{F}_{k_{m},\db_{\ell}}(\Gamma,\lambda) \|^2 \, .
\end{equation}
For the case where $M=1$, and $k_1=k$, we have the single frequency inverse scattering problem given by
\begin{equation}
[\tilde{\Gamma}_k, \tilde{\lambda}_k]=\arg \min_{\Gamma,\lambda} \sum_{\ell=1}^{N_{d}} \| \ub_{k,\db_{\ell}}^{\emph{meas}} -\mathcal{F}_{k,\db_{\ell}}(\Gamma,\lambda)\|^2 \label{eq:single_freq_inv} \, .
\end{equation}

When the scattered field measurements are made at Helmholtz wavenumbers $\leq k$ (for the multifrequency problem $k=k_{M}$), the inverse problems~\eqref{eq:multi_freq_inv}, and~\eqref{eq:single_freq_inv} have the following features. They are inherently ill-posed --- that is, one can expect to stably recover at most $O(k)$ Fourier components of the shape and impedance function due to a version of Heisenberg's uncertainty principle for waves; the stable recovery of these features typically requires measurements at $O(k)$ receivers of the scattered field due to incident waves from $O(k)$ directions; and the objective function becomes increasingly non-convex with increasing $k$, with the size of the local set of convexity in the vicinity of the global minimum shrinking as $O(1/k)$. Moreover, at each iteration in the optimization loop, the evaluation of the objective function requires the solution of $M$ boundary value problems for the Helmholtz equation with $N_{d}$ boundary conditions, where $M$ is the number of frequencies at which the scattered field data is measured and $N_{d}$ is the number of incident plane wave directions. This can be prohibitively expensive, particularly for the multifrequency inverse problem.

In~\cite{chen1995recursive,Chen}, Chen proposed a recursive linearization approach (RLA) to address similar issues that arise in a related multi-frequency inverse problem: that of recovering smoothly varying sound speed profiles of inhomogeneous objects. 
The RLA addresses these issues by reformulating the multifrequency inverse problem as a sequence of constrained single frequency inverse problems and works as a continuation method on the Helmholtz wavenumber, i.e. at wavenumber $k_{j}$, the reconstruction from the previous frequency $k_{j-1}$ serves as an initial guess. The constraints are chosen in order to make the corresponding single frequency problems well-posed. 
By using a continuation method in frequency, the RLA attempts to ensure that sufficiently many features of the unknown sound profile are recovered such that the iterates remain in the local set of convexity of the exact solution for all $k$. This leads to a significant reduction in the number of iterates required to find minimizers of the single frequency inverse problems. The combination of small number iterates required for each optimization problem, and the objective function requiring the solution of only $1$ boundary value problem with $N_{d}$ boundary conditions allows for the efficient reconstruction of the unknown sound speed.
The RLA has been subsequently coupled with fast algorithms for recovering more complicated sound speed profiles in~\cite{Borges2017,Borges2018}, for shape recovery of sound-soft scatters in~\cite{Borges2015}, and for shape recovery of axisymmetric sound-soft scatters in~\cite{Borges2020axis}. 

In this work, we present an extension of this approach for the solution of the multi-frequency inverse problem to recover the shape and impedance, and investigate its behavior for the reconstruction of complicated boundaries for objects greater than $100$ wavelengths in perimeter, and non-smooth impedance functions. The RLA has largely been investigated in the regime where either the number of sensors is $O(k)$ and/or the number of incident directions is $O(k)$. In practice however, measurements are typically made at a fixed number of sensors and a fixed number $O(1)$ of directions of the incident field, and we restrict our attention to this setup.

{\bf Related work:} There are several algorithms for the solution of the inverse problem with classical impedance boundary condition using single frequency data. See, for example, \cite{ivanyshyn2011inverse, kress2001inverse, kress2018inverse, serranho2006hybrid, hettlich1995frechet, lee2007inverse, qin2012inverse, smith1985inverse, yaman2019reconstruction}. Similarly, the inverse problem with generalized impedance boundary conditions has been studied in great detail. See, for example~\cite{aslanyurek2014reconstruction, bourgeois2011stable, bourgeois2012simultaneous, cakoni2012integral, kress2018integral, yang2014reconstruction, kress2019some, guo2015multilayered, yaman2019reconstruction, bourgeois2010identificantion}. For the use of multifrequency data to recover high resolution reconstructions of the shape and sound velocity of a medium, we refer the reader to \cite{bao2005inverse,bao2007inverse,bao2010error, Borges2015, Borges2017, CHAILLAT20124403, chen1995recursive, Chen, sini2012inverse, Farhat2002,harbrecht2007fast,gutman1993regularized}. In particular, a complete review of inverse scattering problems based on multiple frequency data was given in \cite{bao2015inverse}. To the best of our knowledge, we are not aware of previous work on using multifrequency data and recursive linearization to simultaneously find the shape and impedance of an obstacle. 

{\bf Contributions:} The main contribution of this paper is the extension of the RLA to recover the boundary and impedance function of an obstacle with generic boundaries (not necessarily restricted to star-shaped objects) with high resolution,  given multifrequency measurements of its scattered field. This framework can be used for recovering complicated shape and impedance functions and we demonstrate the effectiveness of this approach through several numerical examples. One of the main advantages of this framework is that we do not need to assume any knowledge of the boundary condition. It can be used with scattered data coming from sound-hard obstacles, and with minor modifications for sound-soft obstacles as well. While the RLA tends to obtain high quality reconstructions for most inverse problems, there are still some open questions regarding its convergence behavior which we explore
through additional numerical experiments intended to highlight the mechanism by which the RLA appears to achieve high quality reconstructions when successful. Based on these experiments, we also present a collection of open problems whose resolution could lead to robust and automated tools for the solution of inverse obstacle scattering problems.

%{\color{red}
{\bf Notation:} We present the most common symbols used in this paper in Table \ref{table:symbol}. Whenever we want to refer to a quantity at a wavenumber $k$ and $N_d$ directions, we use the subscript $k$, e.g., $\mathcal{F}_k\coloneqq \left[\mathcal{F}_{k,\db_1}; \cdots;\mathcal{F}_{k,\db_{N_d}}\right]$. Furthermore, ${\boldsymbol u^{\emph{inc}}}$ refers to the vector with coordinates $\left({\boldsymbol u^{\emph{inc}}}\right)_j=u^{\emph{inc}}_{k,\db}(\xb_j)$.%}

\begin{table}
\caption{List of main symbols used in this article.}\label{table:symbol}
{\small
\begin{center}
\begin{tabular}{ll}
\hline
Symbol & Description \\
\hline\hline
$\Omega$                & Closed set representing the impenetrable obstacle\\
$\partial \Omega$     & Boundary of the obstacle $\Omega$\\
$\gamma$                & Parameterization of the curve approximating the boundary $\partial \Omega$ \\
$\Gamma$               & Curve approximating the boundary $\partial \Omega$ \\
$k$                           & Wavenumber (or frequency) of the incident plane wave $u^{\emph{inc}}_{k,\db}$\\
$\db$                        & Incident direction of plane wave $u^{\emph{inc}}_{k,\db}(\xb)$ ($\|\db\|=1$) \\
$\bf r$                       & Receivers locations \\
$u^{\emph{inc}}_{k,\db}$      & Incident plane wave with wavenumber $k$ and incident direction $\db$\\
$u^{\emph{scat}}_{k,\db}$    & Scattered field of the obstacle $\partial \Omega$ generated by $u_{k,\db}^{\emph{inc}}$ \\
${\bf u}^{\emph{meas}}_{k,\db}$  & Vector with coordinates being $u^{\emph{scat}}_{k,\db}$ measured at the receivers\\
$\lambda$                 & impedance boundary function \\
$\mathcal{F}_{k,\db}$  & Forward scattering operator mapping $\partial\Omega$ to $u_{k,\db}^\emph{scat}$ (for given $u_{k,\db}^{\emph{inc}}$) \\
$\partial_\Gamma \mathcal{F}_{k,\db}$ & Frech\'{e}t derivative of $\mathcal{F}_{k,\db}$ with respect to the boundary represented by $\Gamma$ \\
$\partial_\lambda \mathcal{F}_{k,\db}$ & Frech\'{e}t derivative of $\mathcal{F}_{k,\db}$ with respect to the impedance function $\lambda$ \\
$\mathcal{A}_\lambda$  & Space of curves with curvature nearly bandlimited by $\lfloor c_\Gamma k\rfloor$ \\
$\mathcal{A}_\Gamma$ & Space of bandlimited functions with maximum Fourier content  $\lfloor c_\lambda k\rfloor$ \\
$\delta \gamma$       & Update step for the shape in the single frequency inverse problem \\
$\delta \lambda$       & Update step for the impedance function in the single frequency inverse problem \\
$\mathcal{S}_k$           & Single layer potential  at wavenumber $k$\\
$\mathcal{D}_k$           & Double layer potential at wavenumber $k$\\
$\mathcal{K}_k$           & Derivative in the normal direction of the single layer potential at wavenumber $k$\\
$\mathcal{T}_k$           & Derivative in the normal direction of the double layer potential at wavenumber $k$\\
$G^{k}$                     &  Free space Green's function for the 2D Helmholtz equation at wavenumber $k$\\
\hline 
\end{tabular}
\end{center}
}
\end{table}

{\bf Article Outline:} In Section \ref{s:rla}, we present the extension of the RLA for the inverse problem of recovering both shape and impedance. In Section \ref{s:newton}, we present a Gauss-Newton algorithm for the single frequency inverse scattering problem. In Section \ref{s:fwd_problem}, we discuss an integral equation formulation for the solution of the forward scattering problem and its numerical solution. In Section \ref{s:num_res}, numerical examples are presented to illustrate the performance of the method. We present a discussion of known open problems in Section \ref{s:openproblems}. Concluding remarks are made in Section \ref{s:conclusions}.%}

\section{Recursive Linearization approach (RLA)} \label{s:rla}
The recursive linearization approach was introduced in \cite{Chen} for the recovery of sound speed profiles of inhomogeneous objects given far-field/distant measurements of the scattered field at multiple frequencies from multiple incident directions. Without getting into too many details, this problem also requires solving an optimization problem of the form~\ref{eq:multi_freq_inv}, where the unknown function is the sound speed inside the object. 

%which is assumed to be a compactly supported function in $\mathbb{R}^{2}$, denoted by $q$ and the corresponding forward operator $\mathcal{F}_{k_{j},\db_{\ell}}$ requires the solution of a different partial differential equation. 

In the RLA, the multifrequency inverse problem is reformulated as a sequence of constrained single frequency inverse problems, where the initial guess for each single frequency minimization problem is the reconstruction of the solution at the previous frequency. There are three key components to this approach---at each frequency the optimization problems need to be appropriately constrained for stable recovery; using a continuation approach in frequency enables more efficient minimization of increasingly non-convex objective functions;  and finally, at each frequency one just needs to optimize the corresponding single frequency inverse problem which improves the computational efficiency of the method.

\subsection{Constraints on shape and impedance for stable recovery}
First, we set up the constraints on the shape and the impedance of the unknown obstacle in order to obtain a well posed formulation of the multifrequency problem~\cref{eq:multi_freq_inv}. Let $k=k_{M}$ denote the maximum frequency which appears in the objective function. In~\cite{Chen}, Chen observed that for the inverse problem of recovering the sound speed of an inhomogeneous object, one can stably expect to recover $O(k)$ Fourier modes of the sound speed in each direction. This effect is a direct consequence of Heisenberg's uncertainty principle for waves --- sub-wavelength features of the scatterer are present in the evanescent modes of the scattered field and consequently cannot be stably recovered in finite precision arithmetic. Typically, the unknown sound speed of the inhomogeneous object is a smooth compactly supported function in $\mathbb{R}^{2}$, and thus one can directly impose such a constraint on the coefficients of the sine series of the unknown sound speed. 

We would like to impose a similar constraint on the impedance function and the shape of the obstacle for~\cref{eq:multi_freq_inv}. Let $\gamma(t):[0,L] \to \Gamma \subset \bbR^2$ denote an arclength parametrization of a simple closed $C^{3}$ curve $\Gamma$. Let $H(t)$ denote the curvature at $\gamma(t)$ on $\Gamma$. Since $H(t)$ is a periodic function with period $L$, it can be expressed as a Fourier series of the form 
\begin{equation}
H(t) = \sum_{j=-\infty}^{\infty} e^{\frac{2\pi i j t}{L}} \hat{H}_{j} \, .
\end{equation}
One possible approach to restricting the Fourier content of the curve would be through bandlimiting the curvature. However, this requirement turns out to be too restrictive. For example, any star shaped domain of the form $(r(\theta)\cos(\theta), r(\theta) \sin(\theta))$, $0\leq \theta < 2\pi$, which is not a circle does not have bandlimited curvature for any bandlimit. The stringent restrictiveness of bandlimiting the curvature for representing closed curves can in part be explained by the observation made in~\cite{beylkin2014fitting}. Suppose $\Gamma$ is an analytic curve. The Fourier coefficients of the associated curvature decay exponentially. For any $\varepsilon>0$, we can find an appropriate $N$, such that the bandlimited projection of $H(t)$ onto Fourier modes $[-N,N] \subset \mathbb{Z}$, denoted by $H_{N}(t)$ satisfies $\| H_{N} - H \|_{\mathbb{L}^{2}[0,L]} < \varepsilon$. However, the resulting curve with curvature $H_{N}(t)$ need not be closed. This issue can be addressed by adding an $O(\varepsilon)$ correction to the curve using the procedure described in~\cite{beylkin2014fitting}. However, the updated closed curve no longer has bandlimited curvature with the $\mathbb{\ell}^{2}$ norm of the Fourier coefficients in $\mathbb{Z} \setminus [-N,N]$ being $O(\varepsilon)$. With this in mind, we restrict the Fourier content of the curve by constraining $\Gamma \in \mathcal{A}_{\Gamma}(k)$, where the set $\mathcal{A}_{\Gamma}(k)$ is given by
\begin{equation}
\label{eq:cons-gamma}
\mathcal{A}_{\Gamma}(k) = \left\{ \Gamma \, |\, \Gamma \text{ non-intersecting, closed $C^{3}$ curve and } \sqrt{\frac{\sum_{|j|>\floor{c_{\Gamma} k}} |\hat{H}_{j}|^2}{\sum_{j=-\infty}^{\infty} |\hat{H}_{j}|^2}} < \varepsilon_{H}  \right\} \, .
\end{equation}
Here $\hat{H}_{j}$ are the Fourier coefficients of the curvature defined above, and $\varepsilon_{H}$ and $c_{\Gamma}$ are constants.

The task of restricting the Fourier content of the impedance function is much more straightforward. $\lambda(t)$ is also a periodic function with period $L$. Let $\hat{\lambda}_{j}$, $j\in \mathbb{Z}$ denote the Fourier coefficients of $\lambda(t)$,
\begin{equation}
\lambda(t) = \sum_{j=-\infty}^{\infty} e^{\frac{2\pi i j t}{L}} \hat{\lambda}_{j} \, .
\end{equation}
We restrict the Fourier content of $\lambda$ by restricting $\lambda \in \mathcal{A}_{\lambda}(k) \subset \mathbb{L}^{2} [0,L]$, where $\cA_{\lambda}(k)$ is the space of bandlimited functions with maximum Fourier content $\floor{c_{\lambda}k}$, given by
\begin{equation}
\cA_{\lambda}(k) = \{ f \in \mathbb{L}^2[0,L] | \,  , \hat{f}_{j} = 0 \,, \forall |j|>\floor{c_{\lambda}k} \} \, ,
\end{equation}
and $c_{\lambda}$ is a constant.
Given these constraints for the shape of the obstacle and it's impedance, we now seek to solve the constrained multifrequency inverse problem given by
\begin{equation}
\label{eq:multi_freq_inv_cons}
[\tilde{\Gamma},\tilde{\lambda}] = \argmin_{\substack{\Gamma \in \cA_{\Gamma}(k_{M}) \\ 
\lambda \in \cA_{\lambda}(k_{M})}} \sum_{m=1}^{M} \sum_{\ell=1}^{N_{d}} \| \ub_{k_{m},\db_{\ell}}^{\emph{meas}} - \mathcal{F}_{k_{m},\db_{\ell}}(\Gamma,\lambda) \|^2 \, .
\end{equation}

\begin{remark}
The inverse problem of recovering the sound speed of an inhomogeneous object is far more amenable to analysis than corresponding inverse obstacle scattering problems. In particular, in~\cite{Chen}, Chen was able to prove that higher Fourier modes of the unknown sound speed lie in an approximate null space of the forward operator and hence cannot be stably recovered. Due to the non-trivial geometry of the set of closed non-intersecting and approximately bandlimited curves, it is extremely difficult to prove similar results for inverse obstacle scattering problems, and show that~\cref{eq:multi_freq_inv_cons} is well-posed for some $\mathcal{A}_{\Gamma}(k),\mathcal{A}_{\lambda}(k)$.
\end{remark}

\subsection{Continuation in frequency}
\begin{figure}[h!]
\begin{subfigure}[t]{0.33\textwidth}
\center
\includegraphics[width=1\textwidth]{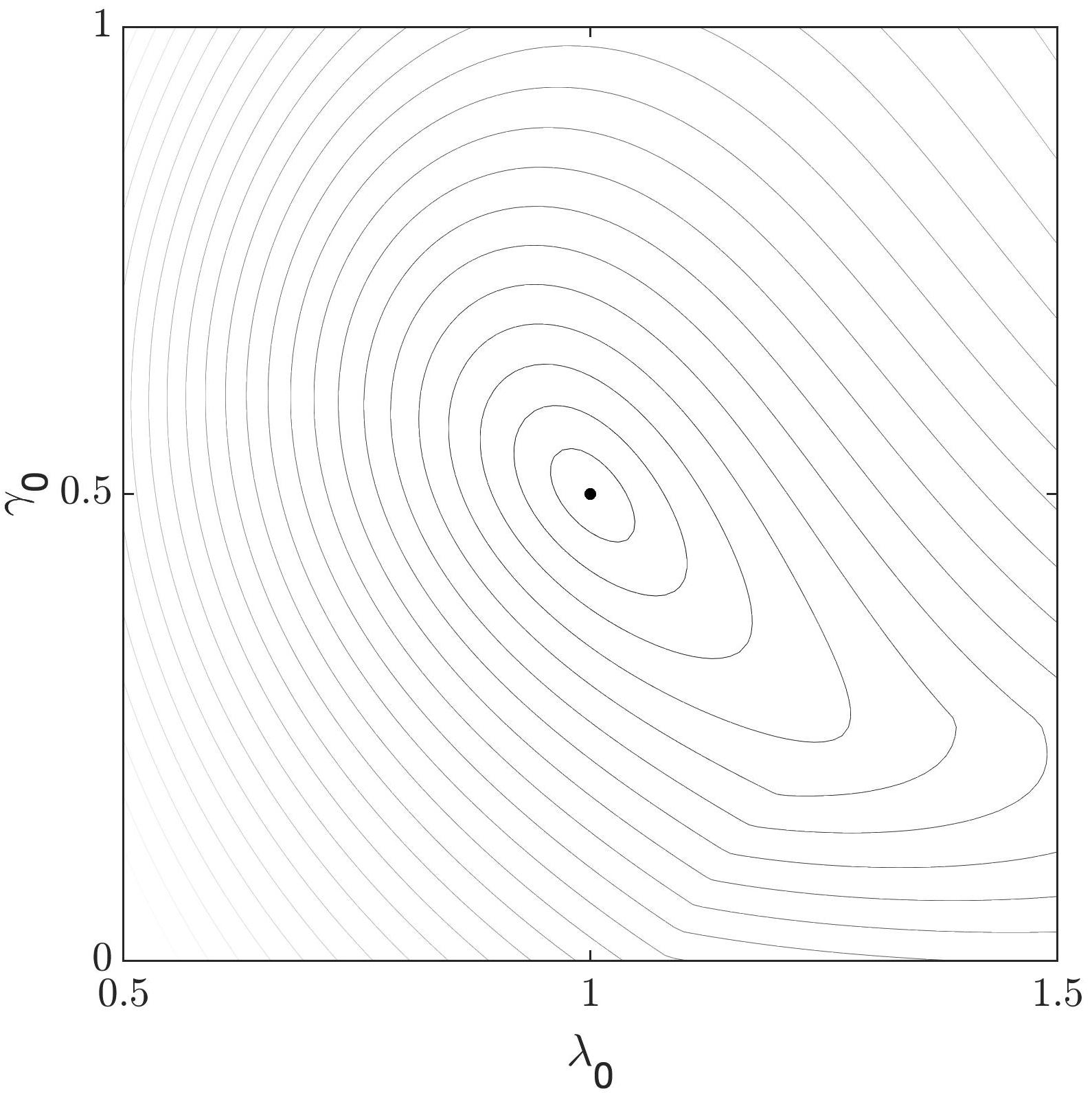}
\caption{$\|{\bf u}_1^{meas}-\mathcal{F}_1(\gamma_0,\lambda_0)\|$}\label{fig:obj_func_k1}
\end{subfigure}
\begin{subfigure}[t]{0.33\textwidth}
\center
\includegraphics[width=1\textwidth]{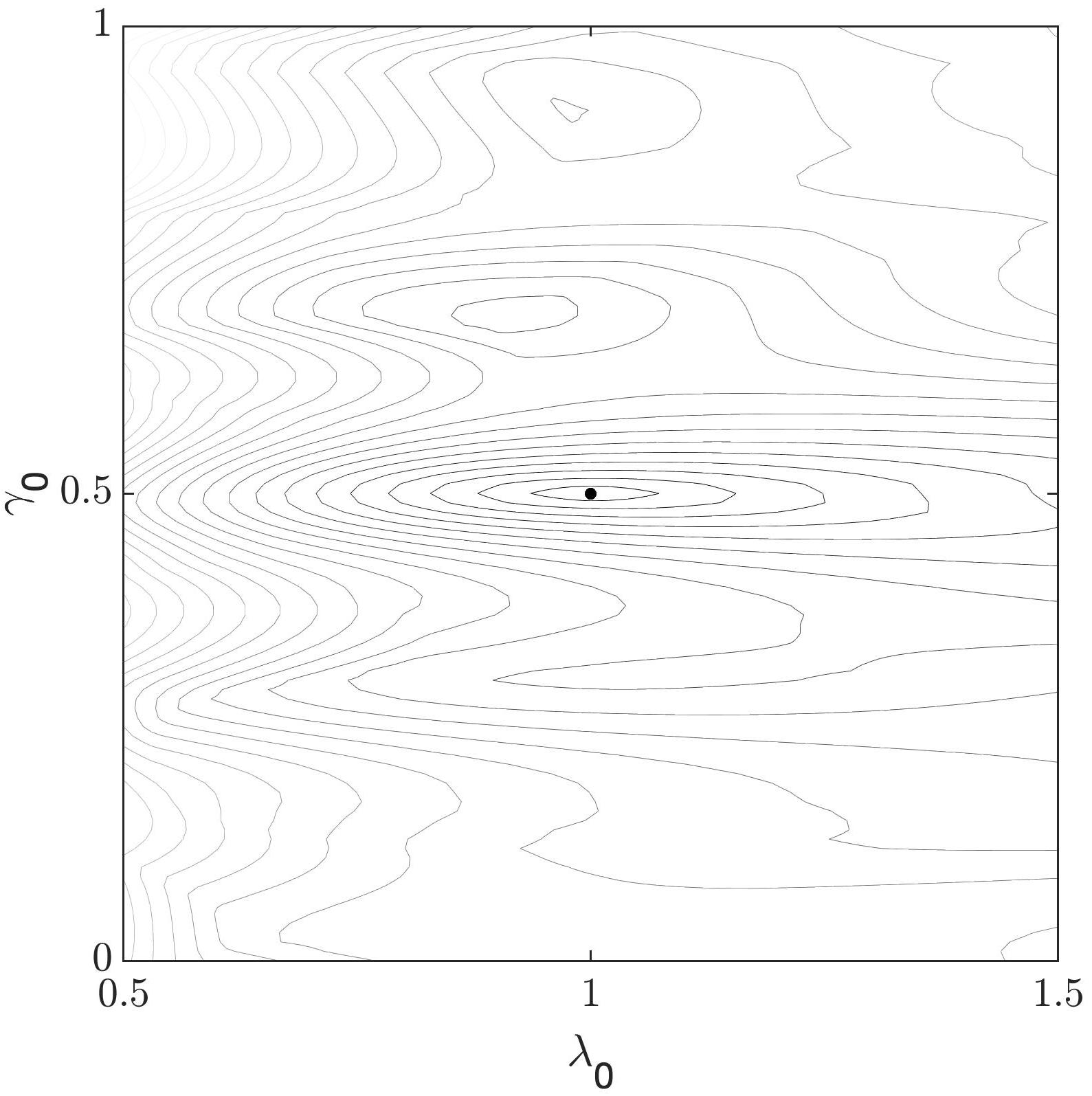}
\caption{$\|{\bf u}_{15}^{meas}-\mathcal{F}_{15}(\gamma_0,\lambda_0)\|$}\label{fig:obj_func_k15}
\end{subfigure}
\begin{subfigure}[t]{0.33\textwidth}
\center
\includegraphics[width=1\textwidth]{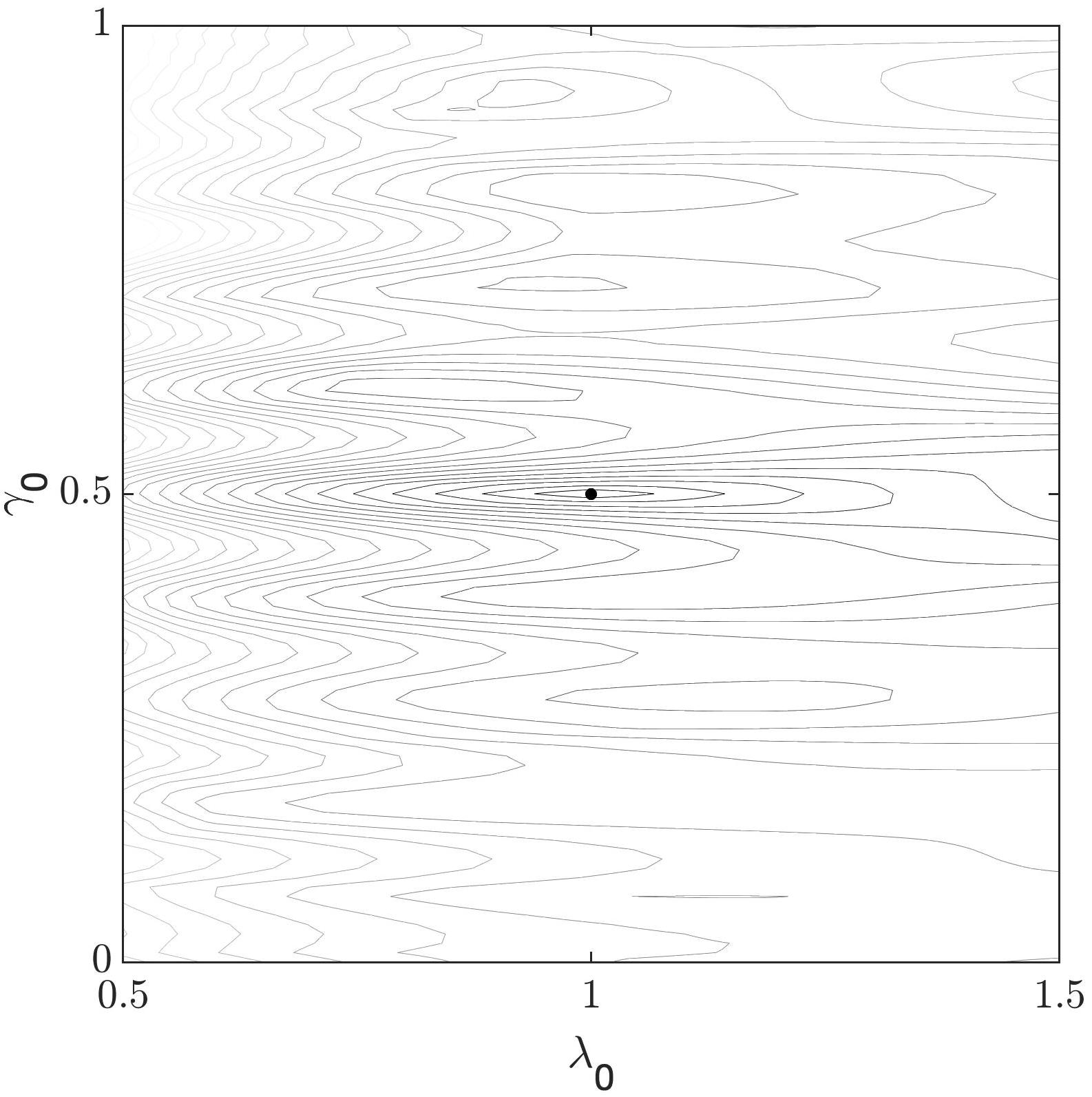}
\caption{$\|{\bf u}_{30}^{meas}-\mathcal{F}_{30}(\gamma_0,\lambda_0)\|$}\label{fig:obj_func_k30}
\end{subfigure}
\caption{Values of the objective function $f_k(\gamma_0,\lambda_0)=\|{\bf u}_k^\emph{meas}-\mathcal{F}_k(\gamma_0,\lambda_0)\|$ for (a) $k=1$, (b) $k=15$ and (c) $k=30$. As the frequency increases, more local minima can be seen near the global minimum of the objective function.}\label{fig:obj_func}
\end{figure}

In order to demonstrate the necessity of using a continuation approach in frequency, we examine the behavior of the objective function in a simple setting. Consider the single frequency objective function $f(\gamma_{0},\lambda_{0}) = \| \ub_{k}^{\emph{meas}} - \cF_{k,\db_{\ell}} (\Gamma,\lambda) \|$, where $\Gamma$ is the boundary of a circle with radius $\gamma_{0}$ and the impedance is a constant denoted by $\lambda_{0}$. For the forward scattering operator $\cF$, we assume that the scattered data is collected from
$N_d=16$ incoming incident waves with incidence direction $\db_j=\pi j/8$, $j=0, \ldots,15$ and measured at $100$ receivers $r_l=10\left(\cos(l\pi/50),\sin(l\pi/50)\right)$, $l=0,\ldots,99$. The measured data $\boldsymbol{u}_{k}^{\emph{meas}}$ is generated for $(\gamma_{0},\lambda_{0}) = (1,0.5)$. In~\cref{fig:obj_func}, we plot the objective function $f$ as we vary $\gamma_{0}$, and $\lambda_{0}$.  As can be seen from the figure, the objective function has several local minima, and spacing between the local minima scales as $O(1/k)$.
The multifrequency objective function behaves similar to the single frequency objective function at the highest frequency.

It follows that the objective function has a large local set of convexity near the solution at low frequencies, while it has a smaller set of convexity near the solution and increasingly many local minima as the frequency is increased. A similar observation was demonstrated for the recovery of three dimension axis-symmetric obstacles in~\cite{Borges2020axis}. Thus, it is essential to have a good initial guess particularly when the measurements are made at high frequencies. One way to obtain a good initial guess is to use multifrequency scattering data and a continuation method in frequency.

A continuation method in $k$ for the constrained inverse problem~\cref{eq:multi_freq_inv_cons} is given by
\begin{equation}
[\tilde{\Gamma},\tilde{\lambda}] = \argmin_{\substack{\Gamma \in \cA_{\Gamma}(k_{M}) \\ 
\lambda \in \cA_{\lambda}(k_{M})}} \sum_{j=1}^{n_{k}} \sum_{\ell=1}^{N_{d}} \| \ub_{k_{j},\db_{\ell}}^{\emph{meas}} - \mathcal{F}_{k_{j},\db_{\ell}}(\Gamma,\lambda) \|^2 = [\tilde{\Gamma}_{M}, \tilde{\lambda}_{M}] \, .
\end{equation}
Here $[\tilde{\Gamma}_{j},\tilde{\lambda}_{j}]$, $j=1,2,\ldots M$ are the solutions of the constrained multifrequency problem up to frequencies $k\leq k_{j}$ given by
\begin{equation}
\label{eq:single_freq_inv_cum}
[\tilde{\Gamma}_{j},\tilde{\lambda}_{j}] = \argmin_{\substack{\Gamma \in \cA_{\Gamma}(k_{j}) \\ 
\lambda \in \cA_{\lambda}(k_{j})}} \sum_{m=1}^{j} \sum_{\ell=1}^{N_{d}} \| \ub_{k_{m},\db_{\ell}}^{\emph{meas}} - \mathcal{F}_{k_{m},\db_{\ell}}(\Gamma,\lambda) \|^2 \, ,
\end{equation}
with initial guess $[\tilde{\Gamma}_{j}^{0}, \tilde{\lambda}_{j}^{0}] = [\tilde{\Gamma}_{j-1}, \tilde{\lambda}_{j-1}]$.

By using a continuation method in frequency, the hope is that the minimizer $[\tilde{\Gamma}_{j-1}, \tilde{\lambda}_{j-1}]$ of the inverse problem with $k\leq k_{j-1}$ lies in the basin of attraction for the inverse problem with $k\leq k_{j}$. Thus, the minimizer for the highly oscillatory loss function associated with the multifrequency problem~\cref{eq:multi_freq_inv_cons}, is constructed in an incremental manner by solving increasingly difficult multifrequency optimization problems, but with increasingly better initial guesses for the minimizers.

\begin{remark}
For the initial guess $[\tilde{\Gamma}_{0}, \tilde{\lambda}_{0}]$ at the lowest frequency $k_{0}$, a very crude approximation of the obstacle and impedance which is close to the location of the true obstacle would be good enough particularly if $k_{0}$ is small.
Alternate strategies for obtaining an initial guess at $k_{0}$ include using multiple signal classification (MUSIC)~\cite{bao2007inverse}, Born approximations~\cite{Colton}, linear sampling methods~\cite{colton_lsm}, and factorization methods~\cite{kirsch_factorization} to name a few.
For the impedance function $\tilde{\lambda}_{0}$, a good initial guess at a low frequencies is a constant function between 0 and 1.
\end{remark}

\subsection{Reduction to single frequency problems}
Recall that evaluating the objective function in~\cref{eq:single_freq_inv_cum} at the highest frequency $k_{M}$, requires the solution of $M$ partial differential equations (PDEs) with $N_{d}$ different boundary conditions for each PDE. Using a large number of frequencies is essential when high resolution features of the unknown obstacles are desired. In this setting,  evaluating the objective function can be prohibitively expensive.
In principle, it is well known that, one could recover both the shape and the impedance of the object from sufficient measurements of the far-field data at a single frequency~\cite{kress2001inverse,Chen}. In particular, at frequency $k$, $O(k)$ measurements of the scattered field from $O(1)$ incident directions are often sufficient for inverse obstacle scattering problems as was demonstrated in~\cite{Borges2015,Borges2020axis}.
Thus, one can improve the computational performance of the method significantly by replacing the multifrequency inverse problem for $k\leq k_{j}$ in~\cref{eq:single_freq_inv_cum}, with the corresponding single frequency inverse problem at $k=k_{j}$ given by
\begin{equation}
\label{eq:single_freq_inv_cons}
[\tilde{\Gamma}_j, \tilde{\lambda}_j]=\argmin_{\substack{\Gamma \in \cA_{\Gamma}(k_{j}) \\ 
\lambda \in \cA_{\lambda}(k_{j})}} \sum_{\ell=1}^{N_{d}} \| \ub_{k_{j},\db_{\ell}}^{\emph{meas}} -\mathcal{F}_{k_{j},\db_{\ell}} \, (\Gamma,\lambda)\|^2 \, .
\end{equation}
The benefit of such a reduction from a computational standpoint is evident since the objective function~\cref{eq:single_freq_inv_cum} requires the solution of $j$ boundary value problems with $N_{d}$ different incident directions, while~\cref{eq:single_freq_inv_cons} requires the solution of just $1$ boundary value problem with $N_{d}$ different incident directions. 

\section{Gauss Newton methods for the single frequency inverse problem~\ref{eq:single_freq_inv_cons}} \label{s:newton}
In the following section, we turn our attention for solving the constrained version of the single frequency inverse problem~\ref{eq:single_freq_inv_cons}.
In an optimization loop, for each guess $[\Gamma, \lambda]$, evaluating the loss function requires the solution of a boundary value problem with $N_{d}$ different boundary conditions. Given the high computation cost for evaluating the objective function, it is imperative to use an optimization method which minimizes the number of function evaluations required. Moreover, using the RLA, we have a good initial guess for the minimizer of~\ref{eq:single_freq_inv_cons} given by the solution of the single frequency problem at the immediately lowest frequency. Finally, derivatives of the objective function can be easily computed at marginal additional cost using the already computed solution operator for the forward problem. In this environment, Gauss-Newton methods are a natural choice for solving the minimization problem. 
  
In the following, we discuss the algorithm for a single incoming wave for notational convenience.
Let $\Gamma^{(j)}, \lambda^{(j)}$, $j=0,1,\ldots$, denote the iterates of the Gauss-Newton algorithm for the boundary $\Gamma$, and impedance function $\lambda$ respectively. As before, let $\gamma^{(j)}(t)$ denote an arc-length parameterization of $\Gamma^{(j)}$. Then the parameterization of the boundary and the impedance function are updated via the formulae
\begin{equation}
\gamma^{(j+1)} = \gamma^{(j)} +\delta \gamma \, , \quad \lambda^{(j+1)} = \lambda^{(j)}+ \delta \lambda \, ,
\end{equation}
where the updates $\delta \gamma$, $\delta \lambda$ satisfy 
\begin{equation} \label{eq:newton_step}
\partial_{\Gamma}\mathcal{F}_{k,\db}(\Gamma^{(j)},\lambda^{(j)}) \delta\gamma + \partial_{\lambda}\mathcal{F}_{k,\db}(\Gamma^{(j)},\lambda^{(j)}) \delta\lambda= {\boldsymbol u}_{k,\db}^\emph{meas} -\mathcal{F}_{k,\db}(\Gamma^{(j)}, \lambda^{(j)}) \, \, .
\end{equation}
Here $\partial_{\Gamma}\mathcal{F}_{k,\db}(\Gamma^{(j)},\lambda^{(j)})$ and $\partial_{\lambda}\mathcal{F}_{k,\db}(\Gamma^{(j)},\lambda^{(j)})$ are the Frech\'et derivatives of the forward operator with respect to the shape $\Gamma$ and impedance function $\lambda$, evaluated at $\Gamma^{(j)}$ and $\lambda^{(j)}$. The precise definition of the Frech\'et derivatives is deferred to~\cref{subsec:frechet} below.  The iterations are repeated until a stopping criteria is reached. The stopping criteria can be the total number of iterations $N_{it}$, the size of the residual $\|{\boldsymbol u}_{k,\db}^\emph{meas}-F_{k,\db}(\Gamma^{(j)},\lambda^{(j)})\|\leq\varepsilon_r$, the size of the update step for the shape $\|\delta\gamma\|_2<\varepsilon_{s,\gamma}$,  the size of the update step for the impedance function $\|\delta\gamma\|_2<\varepsilon_{s,\lambda}$, an increase in the residue $ \|{\boldsymbol u}_{k,\db}^\emph{meas}-F_{k,\db}(\Gamma^{(j)},\lambda^{(j)})\|>\|{\boldsymbol u}_{k,\db}^\emph{meas}-F_{k,\db}(\Gamma^{(j-1)},\lambda^{(j-1)})\|$ or if the updated curve $\Gamma^{(j+1)}$ fails to be in the admissible set.

\subsection{Evaluating the Frech\'et derivatives \label{subsec:frechet}}
The Frech\'et derivatives of $\mathcal{F}_{k,\db}(\Gamma^{(j)},\lambda^{(j)})$ with respect to both the shape $\Gamma^{(j)}$ and the impedance functions $\lambda^{(j)}$ are linear operators. The action of the Frech\'et derivatives on given functions $\delta \gamma$ or $\delta \lambda$, require the solution of the forward impedance
problem albeit with different boundary conditions (the boundary conditions do not arise from point sources or incident plane waves). Explicit expressions for these Frech\'et derivatives have been derived for related forward scattering operators on surfaces in three dimensions, and when the measured data are measurements of the far field pattern of the scattered field~\cite{hettlich1995frechet, Hettlich_1998}  (as opposed to the scattered field measured at receivers). The results extend in a straightforward manner for evaluating the Frech\`et derivatives of $\mathcal{F}_{k,\db}(\Gamma^{(j)},\lambda^{(j)})$. In the following lemma, we state the result for the Frech\`{e}t derivative with respect to the shape $\Gamma$. The proof 
is similar to the results in~\cite{hettlich1995frechet,Hettlich_1998}.

\begin{lemma}\label{thm:der_shape}
Let $\Omega \in \mathbb{R}^2$ be a bounded domain with a $C^{3}$ boundary $\Gamma$, i.e. the associated parameterization has three continuous derivatives. Then the forward problem operator $\mathcal{F}_{k,\db}$ is Fr\'echet differentiable at $\Gamma$, with derivative $\partial_\Gamma \mathcal{F}_{k,\db}(\Gamma,\lambda)\delta\gamma={\boldsymbol v}$, where ${\boldsymbol v}$ is the vector with measurements of the scattered field $v$ at receivers $\boldsymbol{r}_l$ and $v$ satisfies 
\begin{equation}
\begin{cases}
\Delta v+k^2v = 0, \quad &\text{in} ~\Omega, \\
\frac{\partial v}{\partial \nu}+ik \lambda v = k^2 \delta\gamma_\nu u + \frac{d}{ds}\left( \delta\gamma_\nu \frac{du}{ds}\right) -\lambda\delta\gamma_\nu\left(\frac{\partial u}{\partial \nu} -\mathcal{H}u\right) \quad &\text{on} ~\Gamma, \\
\Lim{r\rightarrow \infty} r^{1/2}\left(\frac{\partial v}{\partial r} - i k v\right) = 0, & 
\end{cases}
\end{equation}
where $\delta\gamma_\nu=\delta\gamma\cdot\nu$, $u=u_{k,\db}^\emph{scat}+u^\emph{inc}$, $u_{k,\db}^{\emph{scat}}$ is the solution to equation~\ref{eq:imp_problem}, and $\mathcal{H}$ is the curvature of $\Gamma$.
\end{lemma}

In the following lemma, we state the result for the Frech\`{e}t derivative with respect to the impedance $\lambda$. The proof 
is similar to the results in~\cite{hettlich1995frechet,Hettlich_1998}.

\begin{lemma}\label{thm:der_imp}
Let $\Omega \in \mathbb{R}^2$ be a bounded domain with a $C^{3}$ boundary $\Gamma$. Then the forward problem operator $\mathcal{F}_{k,\db}$ is Fr\'echet differentiable at $\Gamma$, with derivative $\partial_\lambda \mathcal{F}_{k,\db}(\Gamma,\lambda)\delta\lambda={\boldsymbol w}$, where ${\boldsymbol w}$ is the vector with measurements of the scattered field $w$ at receivers $\boldsymbol{r}_l$ and $w$ satisfies 
\begin{equation}
\begin{cases}
\Delta w+k^2w = 0, \quad &\text{in} ~\Omega, \\
\frac{\partial w}{\partial \nu}+ik \lambda w = -ik\delta \lambda u \quad &\text{on} ~\Gamma, \\
\Lim{r\rightarrow \infty} r^{1/2}\left(\frac{\partial w}{\partial r} - i k w\right) = 0 \, . & 
\end{cases}
\end{equation}
\end{lemma}

It follows that the marginal cost of computing the action of Frech\`{e}t derivatives on given perturbations is comparable to the cost of evaluating the objective function, and this cost is significantly smaller if the solution operator is precomputed or compressed. 

\subsection{Computing the updates $\delta \gamma$ and $\delta \lambda$}
In order to compute the updates $\delta \gamma$, and $\delta \lambda$ we need to solve equation~\ref{eq:newton_step} in an appropriate basis for $\delta \gamma$ and $\delta \lambda$. In the remainder of the manuscript, for ease of notation, let $N_{\Gamma} = \floor{c_{\Gamma}k}$, and $N_{\lambda} = \floor{c_{\lambda} k}$. Since $\delta \lambda \in \mathcal{A}_{\lambda}(k)$, it is natural to express $\delta \lambda$ through it's Fourier representation as
\begin{equation*}
\delta\lambda(t)= \delta \lambda_{0}^{(c)} + \sum_{\ell=1}^{N_{\lambda}} \left[ \delta \lambda_{\ell}^{(c)}\cos\left(\frac{2\pi \ell t}{L} \right) + \delta \lambda_{\ell}^{(s)} \sin\left(\frac{2\pi \ell t}{L}\right) \right] \, ,
\end{equation*}
where $L$ is the length of the curve $\Gamma^{(j)}$, and $\delta \lambda_{0}^{(c)}$, and $\delta \lambda_{\ell}^{(c)}, \delta \lambda_{\ell}^{(s)}$, $\ell=1,2\ldots N_{\lambda}$ are real constants.

The update associated with the boundary requires a little more care. The component $\delta \gamma \cdot \tau$ does not change the shape of the obstacle and only changes the parameterization of the curve. Moreover, as is evident from the formula of the Frech\`et derivative $\partial_{\Gamma} \cF_{k,\db}$, any tangential update to the curve lies in the null space of $\partial_{\Gamma} \cF_{k, \db}$. Thus, we parametrize the update of the shape as $\delta \gamma = h(t) \bnu^{(j)}(t)$, where $\bnu^{(j)}(t)$ is the normal to the curve at $\gamma^{(j)}(t)$, and $h(t)$ is a scalar function. However, unlike the impedance update, there doesn't exist a choice of $h(t)$ which would guarantee that the curvature of the updated curve $\gamma^{(j+1)} = \gamma^{(j)} + h(t) \bnu^{j}(t)$ would be constrained to $\cA_{\Gamma}(k)$. This is partly due to the curvature being a non-linear function of the update $h(t)$, and also for generic geometries which are not star shaped, it is difficult to guarantee that the updated curve would not be self-intersecting. Even if the curvature of $\gamma^{(j)}$ is nearly $0$, $h(t)$ must still be restricted to bandlimited functions on the curve with bandlimit $N_{\Gamma}$ in order for $\gamma^{(j+1)}(t) \in \cA_{\Gamma}(k)$. Thus we assume that the update 
$h$ takes the form
\begin{equation}
h(t)=\delta\gamma^{(c)}_0+\sum_{\ell=1}^{N_\gamma} \left[\delta\gamma^{(c)}_\ell \cos\left(\frac{2\pi \ell t}{L}\right)+\delta\gamma^{(s)}_\ell \sin\left(\frac{2\pi \ell t}{L}\right)\right],
\end{equation}
where $\delta\gamma^{(c)}_{0}$, and $\delta\gamma^{(c)}_{j}$, $\delta\gamma^{(s)}_{j}$, $\ell=1,2,\ldots N_{\Gamma}$ are real constants.

For each basis element of $\delta \lambda(t) = \cos{(2\pi \ell t/L)}$, $\sin{(2\pi \ell t/L)}$, and of $\delta \gamma(t) = \cos{(2\pi \ell t/L)} \bnu^{(j)}(t)$,  
$\sin{(2\pi \ell t/L)} \bnu^{(j)}(t)$, we compute the action of the respective Frech\`{e}t derivative to form a discrete linear system.
The unknown coefficients $\delta\lambda_{0}^{(c)}$, $\delta \lambda_{\ell}^{(c,s)}$, $\ell=1,2,\ldots N_{\lambda}$,   $\delta \gamma^{(c)}_{0}$, and $\delta \gamma_{\ell}^{(c,s)}$, $\ell=1,2,\ldots N_{\Gamma}$ are then obtained by solving the discrete linear system in a least square sense.

However, there is no guarantee that the updated curve satisfies the constraint $\Gamma^{(j+1)}(t) \in \cA_{\Gamma}(k)$. 
If $h(t)$ results in $\Gamma^{(j+1)}(t) \not \in \cA_{N_{\gamma}, \varepsilon}$, we filter the update coefficients $\delta \gamma_{\ell}^{(c)}$, and $\delta \gamma_{\ell}^{(s)}$ using a Gaussian filter with variance $\sigma$ as follows:
\begin{equation}
\delta \gamma_{\ell}^{(c)} \to \delta \gamma_{\ell}^{(c)} e^{-\frac{\ell^2}{N_{\gamma}^2 \sigma^2}} \, , \quad \delta \gamma_{\ell}^{(s)} \to \delta \gamma_{\ell}^{(s)} e^{-\frac{\ell^2}{N_{\gamma}^2 \sigma^2}} \, ,\quad \ell =1,2,\ldots N_{\gamma} \, .
\end{equation}
We repeatedly apply the filter with smaller values of $\sigma$ until we obtain an updated curve $\Gamma^{(j+1)} \in \cA_{N_{\Gamma}, \varepsilon}$.  

\begin{remark}
From the perspective of the Gauss-Newton method, the filtering procedure results in a change in the update direction used for computing $(\Gamma^{(j+1)}, \lambda^{(j+1)})$. This could lead to potential issues for the optimization method. One could remedy this issue by computing an additional update to the impedance while holding the shape of the boundary fixed, if the boundary update is filtered. However, in practice we observe that even though this additional step reduces the residual, there is a marginal reduction in the quality of the reconstruction. and hence we do not include the correction step in our algorithm.
\end{remark}

We summarize the Gauss-Newton algorithm for the constrained single frequency minimization problem below:
\begin{algorithm}
\caption{Gauss-Newton method for shape and impedance}
\label{alg:GN}
\begin{algorithmic}[1]
\STATE{{\bf Input:} Scattered field measurements ${\ub}^\emph{meas}_k$, initial guesses $\gamma^{(0)}$ and $\lambda^{(0)}$, parameters $\alpha$, $\sigma_\gamma$, $\sigma_\lambda$, $N_{it}$, $\epsilon_r$, and $\epsilon_s$.}
\STATE{Set $j=0$, $\gamma=\gamma^{(0)}$, $\lambda=\lambda^{(0)}$, $\|\delta\gamma\|=2 \epsilon_s$ and $\|\delta\lambda\|=2 \epsilon_s$, $r_{0} = \min{(2\pi/k,1)}$}
\WHILE{$j<N_{it}$ {\bf and}  $\|{\bf u}^\emph{meas}_k- \mathcal{F}_k(\gamma)\|> \epsilon_r$ {\bf and} $\|\delta\gamma\|> \epsilon_s$ {\bf and} $\|\delta\lambda\|> \epsilon_s$}
\STATE{Calculate $\mathcal{F}_k(\Gamma^{(j)},\lambda^{(j)})$ and the derivatives $\partial_\Gamma \mathcal{F}_k(\Gamma^{(j)},\lambda^{(j)})$ and $\partial_\lambda \mathcal{F}_k(\Gamma^{(j)},\lambda^{(j)})$.}
\STATE{Solve $\partial_{\Gamma}\mathcal{F}_k(\Gamma^{(j)},\lambda^{(j)}) \delta\gamma + \partial_{\lambda}\mathcal{F}_k(\Gamma^{(j)},\lambda^{(j)}) \delta\lambda = {\bf u}_k^\emph{meas} -\mathcal{F}_k(\gamma^{(j)},\lambda^{(j)})$.}
\STATE{$\sigma=1$}
\STATE{$\gamma^{(j+1)}\leftarrow\gamma^{(j)}+ \delta\gamma$}
\STATE{$\lambda^{(j+1)}\leftarrow\gamma^{(j)}+ \delta\lambda$}
\WHILE{$\Gamma_{j+1} \not \in \mathcal{A}_{\Gamma}(k)$}
\STATE{$\widetilde{\delta\gamma}\leftarrow$ filter($\delta\gamma$,$\sigma$)}
\STATE{$\sigma \leftarrow$ $\sigma/10$}
\STATE{$\gamma^{(j+1)}\leftarrow\gamma^{(j)}+\widetilde{\delta\gamma}$}	
\ENDWHILE
\STATE{$j\leftarrow j+1$}
\ENDWHILE
\end{algorithmic}
\end{algorithm}

\section{Forward Scattering Problem}\label{s:fwd_problem}
At each Newton iteration, we need to solve the impedance boundary value problem~\ref{eq:imp_problem} for a given boundary $\Gamma$, and impedance function $\lambda$ to evaluate the objective function $\mathcal{F}_{k,\db}(\Gamma^{(j)},\lambda^{(j)})$ and its Frech\`et derivatives. Note that the boundary data need not necessarily be of the form $$-\left(\frac{\partial \uinc}{\partial \nu} + ik \lambda \uinc \right)$$ and could be a generic smooth function. In this work, we use integral equation methods to solve this boundary value problem. 

Let $G_{k}(\bx,\by)$ denote the Green's function of the Helmholtz equation with wave number $k$ given by
\begin{equation}
G_{k}(\bx,\by) = H_{0}^{(1)}(k|\bx-\by|) \, ,
\end{equation}
where $H_{0}^{(1)}(z)$ is the Hankel function of the first kind of order $0$. Let $\cS_{k}[\sigma](\bx)$, and $\cD_{k}[\sigma](\bx)$ denote the corresponding single and double layer potentials  given by
\begin{equation}
\mathcal{S}_{k}[\sigma](\bx) = \int_{\Gamma} G_{k}(\bx,\by) \sigma(\by) ds(\by) \, , \quad \text{and} \quad \mathcal{D}_{k}[\sigma] = 
\int_{\Gamma} \nabla G_{k}(\bx,\by) \cdot \boldsymbol{\nu}(\by) \sigma(\by) ds (\by) \, .
\end{equation}
We use the following regularized combined field integral equation representation for the potential $u$
\begin{equation}
u = \cS_{k} [\sigma] + i k \cD_{k} \cS_{i|k|} [\sigma] \, ,
\end{equation}
where $\sigma$ is an unknown density.
This representation was first used for sound-hard scatterers in~\cite{bruno_2012}.
By construction, $u$ satisfies the Helmholtz equation in $\mathbb{R}^{2} \setminus \Omega$, and the Sommerfeld radiation condition. Applying the boundary condition $-(\partial u/\partial n + ik\lambda u) =f$ along $\Gamma$, we get the following integral equation for the unknown density $\sigma$, 

\begin{equation}
\label{eq:inteq0}
\left( -\frac{1}{2} \mathcal{I}+ \mathcal{S'}_{k}^{\textrm{PV}} + ik \mathcal{T} \cS_{ik}  + ik  \lambda \left( \cS_{k} + ik\cD^{\textrm{PV}} \cS_{i|k|} +\frac{ik}{2}\cS_{i|k|} \right) \right)\sigma= f \, ,\quad \bx \in \Gamma \, .
\end{equation}
Here $\mathcal{S}_{k}'^{\textrm{PV}}$ is the principal value of the normal derivative of the single layer potential given by
\begin{equation}
\mathcal{S}_{k}'^{\textrm{PV}}[\sigma] = \textrm{ p.v.} \left(\nu (\bx) \cdot \int_{\Gamma} \nabla G_{k}(\bx,\by)  \sigma(\by) ds (\by)  \right)\, ,
\end{equation}
$\mathcal{D}_{k}^{\textrm{PV}}$ is the principal value of the double layer potential given by
\begin{equation}
\mathcal{D}_{k}^{\textrm{PV}}[\sigma] = 
\textrm{p.v.} \int_{\Gamma} \nabla G_{k}(\bx,\by) \cdot \nu(\by) \sigma(\by) ds (\by) \, ,
\end{equation}
and $\mathcal{T}_{k}[\sigma]$ is the finite part of the normal derivative of the double layer potential given by
\begin{equation}
\mathcal{T}_{k}[\sigma] = 
\textrm{f.p.} \left( \nu(\bx) \cdot \int_{\Gamma} \nabla G_{k}(\bx,\by) \cdot \nu(\by) \sigma(\by) ds (\by) \right) \, .
\end{equation}
Using Calderon identities, equation~\ref{eq:inteq0} can be further simplified to the following second kind integral equation
\begin{equation}\label{eq:inteq}
\begin{aligned}
\Bigg[ -\frac{2 + ik}{4} \mathcal{I}+ \mathcal{S'}_{k}^{\textrm{PV}} &+ ik \left( \left( \mathcal{T}_{k} - \mathcal{T}_{i|k|} \right)  \cS_{i|k|}  + \left(\mathcal{S'}_{i|k|}^{\textrm{PV}}\right)^2 \right) +  \\
&ik  \lambda \left( \cS_{k} + ik \cD^{\textrm{PV}} \cS_{i|k|} +\frac{ik}{2}\cS_{i|k|} \right) \Bigg]\sigma= f \, ,\quad \bx \in \Gamma \, ,
\end{aligned}
\end{equation}
This integral equation can be viewed as the right preconditioned analog of the corresponding left-preconditioned equation presented in~\cite{colton2013integral}, with the additional observation of~\cite{bruno_2012} that it is advantageous to use $\cS_{i|k|}$ as a preconditioner as opposed to $\cS_{0}$. The proofs discussed in~\cite{colton2013integral,bruno_2012} can be extended to show that the integral equation~\ref{eq:inteq} is invertible as long as $k,\lambda$ are real, and $\lambda>0$.

\subsection{Discretization details and numerical solution of the integral equation}
The updated curve after every Newton iteration will not have an arclength parametrization and over time, this might result in a significant deterioration in the quality of the parametrization. In order to avoid this situation, we reparameterize the boundary and construct an arc-length parameterization at every iteration. Moreover, since the curvature of the boundary $\Gamma$ is nearly bandlimited, and the impedance is bandlimited, where the constants $c_{\Gamma},c_{\lambda}$ are typically being chosen $O(1)$, it suffices to sample the geometry at $40$ points per wavelength at each frequency to get a fairly accurate solution. We use an equispaced discretization for representing the curve $\Gamma$, and use Gauss-trapezoidal rule of order 16 proposed by Alpert \cite{alpert} for discretizing the operators $\cS$, $\cD, \cS'$, and $(\mathcal{T} -\mathcal{T}_{i|k|})$. Once $\sigma$ is available, we first evaluate $\tilde{\sigma} = \cS_{i|k|}[\sigma]$ using the appropriate discretized matrix, and then the scattered field at the receivers can be evaluated by discretizing the single and double layer potentials (acting on $\sigma$ and $\tilde{\sigma}$ respectively) using the trapezoidal rule since the targets are not close to the boundary.

The same PDE needs to be solved for $N_{d} \times (N_{\Gamma} + N_{\lambda} + 1)$ different boundary conditions at each Newton iteration. In this work, for ease of implementation, we form the dense matrices and invert it using an $LU$ factorization. The computation cost scales like $O(N^3)$ where $N$ is the number of discretization points. However, one could apply fast direct solvers to both form and apply the inverse in $O(N)$ time~\cite{chandrasekaran2006fast, chandrasekaran2006fast1, gillman2012direct, greengard2009fast, ho2012fast, martinsson2005fast,bebendorf2005hierarchical, bormhierarchical, borm2003introduction} to accelerate the problem at higher frequencies if needed. 

\section{Numerical Experiments}\label{s:num_res}

We present three numerical examples to illustrate the performance of the RLA for the recovery of shape and impedance given measurements of the scattered field.  The performance of the RLA is sensitive to the choice of $N_{\Gamma} = \floor{c_{\Gamma}k}$, and $N_{\lambda} = \floor{c_{\lambda}k}$ which define the subsets in which we solve the single frequency minimization problems. Setting $c_{\Gamma}$, $c_{\lambda}$ too small can result in the algorithm not ending up in the appropriate basin of attraction as we march in frequency and converging to a local minimum. On the other hand, setting them to be too large can result in instabilities in the Gauss-Newton iteration due to the inherent ill-posedness of the problem. Moreover, the amount of information that can be stably extracted from the data at a given frequency is fixed, and needs to be appropriately distributed between recovering modes of the shape and the impedance. 

In~\cref{sec:ex1}, we present heuristics for choosing $c_{\Gamma}, c_{\lambda}$ which address the concern above. This is done by studying the performance of the algorithm in three cases: reconstructing the obstacle assuming the impedance function is known; reconstructing the impedance assuming the shape of the obstacle is known; and finally reconstructing both the shape of the obstacle and the impedance function. Through this experiment, we also study the behavior of the algorithm when the data is noisy, but still in the high signal to noise regime.

A remarkable and surprising feature of the approach is its ability to recover the shape of the obstacle for both sound-soft, and sound-hard scatterers. This feature obviates the need for knowing appropriate boundary conditions for the object being reconstructed and is demonstrated in~\cref{sec:ex2}. Finally, in~\cref{sec:ex3}, we demonstrate the ability of the RLA to recover high frequency features of the shape and the impedance function. A list of figures associated with these experiments is presented in~\cref{tab:examples}. 
%\begin{table}[!htbp]	
%	\begin{center}
%		\begin{tabular}{|c | l | c|}\hline
%			Example &  Description & Figures \\ \hline\hline
%			\multirow{ 2}{*}{1} & Recovering the a) shape, b) the impedance function and  & \ref{fig:ex1a_final} ,\ref{fig:ex1b_final}, \ref{fig:ex1c_final_shape}, \ref{fig:ex1c_final_impedance}, \\
%			   & c) both shape and impedance function & \ref{fig:ex1c_final_residues}, \ref{fig:ex1c_final_best_recons}\\\hline			
%			2 & Using Neumann and Dirichlet scattered data & \\\hline
%			3 & Recovering high frequency information &  \\\hline
%		\end{tabular}
%	\end{center}
%	\caption{List of numerical examples with respective results.}\label{tab:examples}
%\end{table}

\begin{table}[!htbp]	
	\caption{List of figures with description of results and related section in the article.}\label{tab:examples}
	\begin{center}
		\begin{tabular}{|c | l | c|}\hline
			Example &  Description & Figures \\ \hline\hline
			\multirow{ 2}{*}{5.1} & Recovering the: shape, the impedance function and  & \multirow{ 2}{*}{3--8} \\
			   & both shape and impedance function & \\\hline
			5.2 & Using Neumann and Dirichlet scattered data & 9--10 \\\hline
			5.3 & Recovering high frequency information & 11--18 \\\hline
		\end{tabular}
	\end{center}
\end{table}

For each example, unless stated otherwise, we assume that scattered field measurements are made for $M$ frequencies, $k_\ell=k_0+(\ell -1)\delta k $, $\ell=1,\ldots M$, with $\delta k = 0.25$, and $k_{0}=1$. Let $\kmax = 1+(M-1)\delta k$ denote the maximum frequency for which the data is available. At each frequency, the data is generated by $N_d=16$ incident waves, with incidence directions $\db_j=2j\pi/N_d$, $j=1,\ldots,N_d$. The data is measured at $N_r=100$ receivers located at points $\br_m=10(\cos(2\pi m/N_r),\sin(2\pi m/N_r))$, for $m=1,\ldots,N_r$. The parameters associated with the stopping criterion for Gauss-Newton algorithm (~\cref{alg:GN}) are set as follows: maximum number of iterations $N_{it} = 200$; residual tolerance $\varepsilon_{r} = 10^{-3}$; and tolerance of the size of the update step for the impedance function $\varepsilon_{s,\lambda} = 10^{-3}$. We do not impose any restriction on the update step for the boundary. For defining the constraint set for curves $\cA_{\Gamma}(k)$, we set the maximum allowed $L^{2}$ energy of the high Fourier components of the curvature denoted by $\varepsilon_{H}$ in~\cref{eq:cons-gamma} to $10^{-3}$. For $k=1$, we initialize the domain to be a unit circle centered at the origin, and the impedance to $\lambda = 1$. In a slight abuse of notation, let $\varepsilon_{r}$ also denote the relative residual 
\begin{equation}
\varepsilon_{r} = \frac{\left\| \ub_{k}^{\emph{meas}} - \mathcal{F}_{k}(\tilde{\Gamma}_{k}, \tilde{\lambda}_{k}) \right\|}{\left\| \ub_{k}^{\emph{meas}} \right\|} \, ,
\end{equation}
where $\tilde{\Gamma}_{k},\tilde{\lambda}_{k}$ are the reconstructions for $\Gamma,\lambda$, when the Gauss Newton algorithm for the single frequency minimization problem with frequency $k$ is terminated. Let $\varepsilon_{\lambda}$ denote the error in the reconstructed $\lambda$. Note that even though the reconstruction is not defined on the underlying obstacle, we can rescale the two functions to be defined on $(0,2\pi)$ and then compute the error on the interval $(0,2\pi)$. More precisely, let $L$ denote the length of the boundary of the obstacle $\pa \Omega$, and let $\tilde{L}$ denote the length of the curve $\tilde{\Gamma}$, then
\begin{equation}
\varepsilon_{\lambda} = \sqrt{\left(\int_{0}^{2\pi} \left| \lambda\left( \frac{Lt}{2\pi} \right) - \tilde{\lambda}\left(\frac{\tilde{L} t}{2\pi}\right) \right|^2 \, dt\right)}
\end{equation}

In order to avoid inverse crimes, we solve the forward problem using a different number of points than the inverse problem. To generate the scattered field measurements, we solve the problem by discretizing the boundary using $50$ points per wavelength, while using $40$ points per wavelength for the reconstruction. In practice, we observe that solving the forward problem using a different representation (for example, $u = \cS[\sigma] + ik \cD[\sigma]$) tends to have no impact on the reconstruction obtained.

\subsection{Selecting $c_{\Gamma}$, $c_{\lambda}$  \label{sec:ex1}} 
In order to obtain heuristics for selecting $c_{\Gamma}$, and $c_{\lambda}$, we appeal to the existing heuristics available for the two limiting problems: recovering the shape of the obstacle when the impedance is known, and recovering the impedance of the obstacle when the shape is known. The former is similar to using the RLA for inverse obstacle scattering problem for sound-soft scatterers discussed in~\cite{Borges2015}, while the latter is similar to recovering the sound speed of an inhomogeneous object discussed in~\cite{Chen}. Both of those heuristics indicate that setting $c_{\Gamma} = c_{\lambda} = 2$ should result in good reconstructions. We test this hypothesis for a star shaped obstacle whose boundary is parametrized 
as  $\gamma(t)=r(\theta)(\cos(\theta),\sin(\theta))$, $0\leq t<2\pi$, with 
\begin{equation*}
r(\theta)=1+0.2\cos(3\theta)+0.02\cos(4\theta)+0.1\cos(6\theta)+0.1\cos(8\theta) \, ,
\end{equation*}
 and the impedance function $\lambda:[0,2\pi]\rightarrow\mathbb{R}$ is given by the trigonometric function
 \begin{equation*}
 \lambda(t)= 1 + 0.1 \cos(t) + 0.02 \cos(9t) \, .
 \end{equation*}
 For all the experiments in this section, we add $2\%$ noise to the scattered field measurements and set
 \begin{equation}
 \label{eq:noise}
 u^\emph{scat}=u^\emph{scat}+ 0.02 |u^\emph{scat}|\frac{\Phi}{|\Phi|} \, ,
 \end{equation}
 where $\Phi = \phi_{1} + i \phi_{2}$, and $\phi_{1},\phi_{2}$, are i.i.d Gaussians with mean $0$ and standard deviation $1$. 
 
\begin{figure}[h]
 \center
\includegraphics[width=\linewidth]{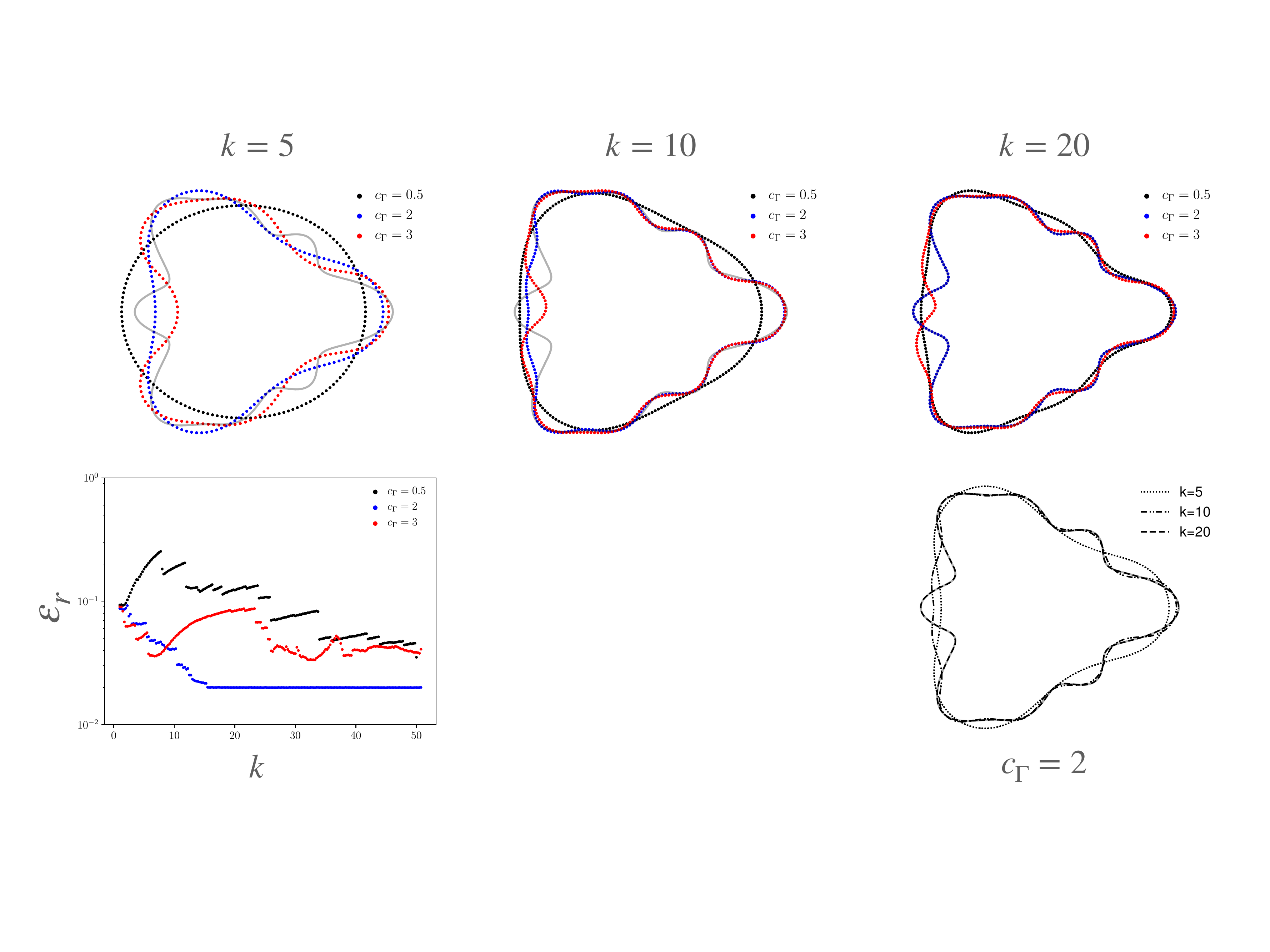}
\caption{{\bf (Section 5.1)} Results for the shape reconstruction at $k=5$, $10$ and $20$ using $c_\Gamma=0.5$, $2$ and $3$ when recovering only the shape of the obstacle. In this case, the impedance function is known.}\label{fig:ex1a_final}
\end{figure}

First we reconstruct the shape of the obstacle assuming the impedance to be known using~\cref{alg:GN} with $c_{\Gamma} =2$, $c_{\Gamma}=3$, and $c_{\Gamma} = 0.5$, and $\kmax=50$. In this example, we impose an additional stopping criterion on the step size for the shape update, $\varepsilon_{s,\gamma} = 10^{-3}$.
In~\cref{fig:ex1a_final}, we compare the reconstructions of the boundary of the obstacle for these configurations at $k=5$, $k=10$, and $k=20$, plot the relative residue $\varepsilon_{r}$, and plot the reconstruction of the shape as a function of frequency for $c_{\Gamma}=2$. Similar to the performance of the RLA for sound-soft obstacles in~\cite{Borges2015}, we observe that the method results in high quality reconstructions for $c_{\Gamma}=0.5,2$, and the method results in instabilities in the reconstruction for $c_{\Gamma}=3$. As expected, the rate of convergence of the reconstruction is much slower for $c_{\Gamma}=0.5$ as compared to $c_{\Gamma}=2$.

\begin{figure}[h]
 \center
\includegraphics[width=\linewidth]{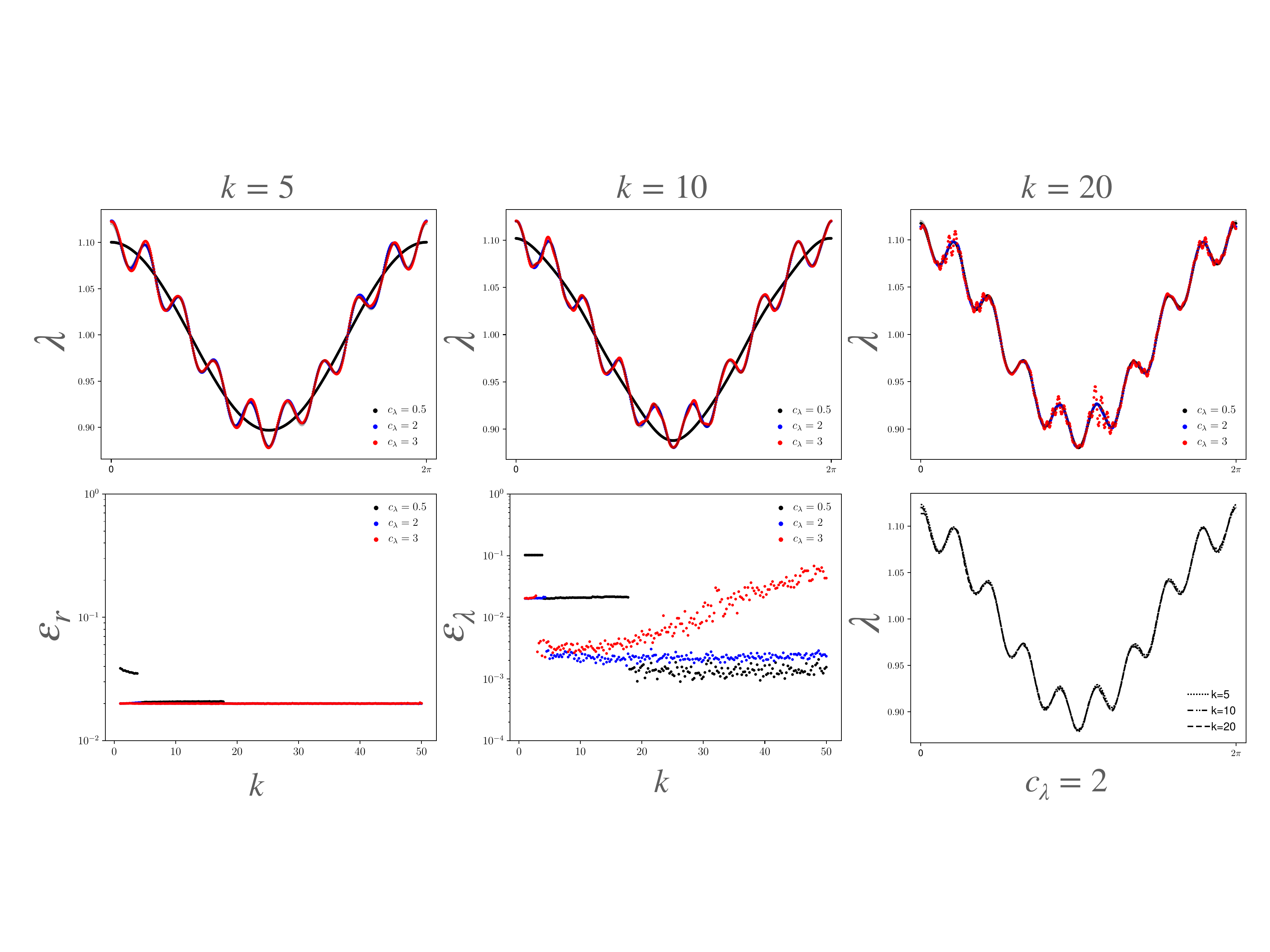}
\caption{{\bf (Section 5.1)} Results for the impedance function reconstruction at $k=5$, $10$ and $20$ using $c_\lambda=0.5$, $2$ and $3$ when recovering only the impedance function. In this case, the shape of the obstacle is known.}\label{fig:ex1b_final}
\end{figure}

Next, we reconstruct the impedance of the obstacle assuming that the shape of the obstacle is known using~\cref{alg:GN} with $c_{\lambda} = 0.5,2$, and $3$, and $\kmax=50$. This problem is significantly easier as compared to the inverse problems involving the shape of the obstacle due to the simplicity of the function space and constraint set for $\lambda$. In~\cref{fig:ex1b_final}, we compare the reconstructions for these configurations at $k=5$, $k=10$, and $k=20$, plot the relative residue $\varepsilon_{r}$, the error in impedance $\varepsilon_{\lambda}$, and plot the reconstruction of the impedance as a function of frequency for $c_{\lambda}=2$.  Similar to the reconstruction of the shape, we observe instabilities for $c_{\lambda}=3$, but are able to recover the impedance function up to the level of noise in the data with $c_{\lambda}=0.5,2$.

\begin{figure}[h]
 \center
\includegraphics[width=\linewidth]{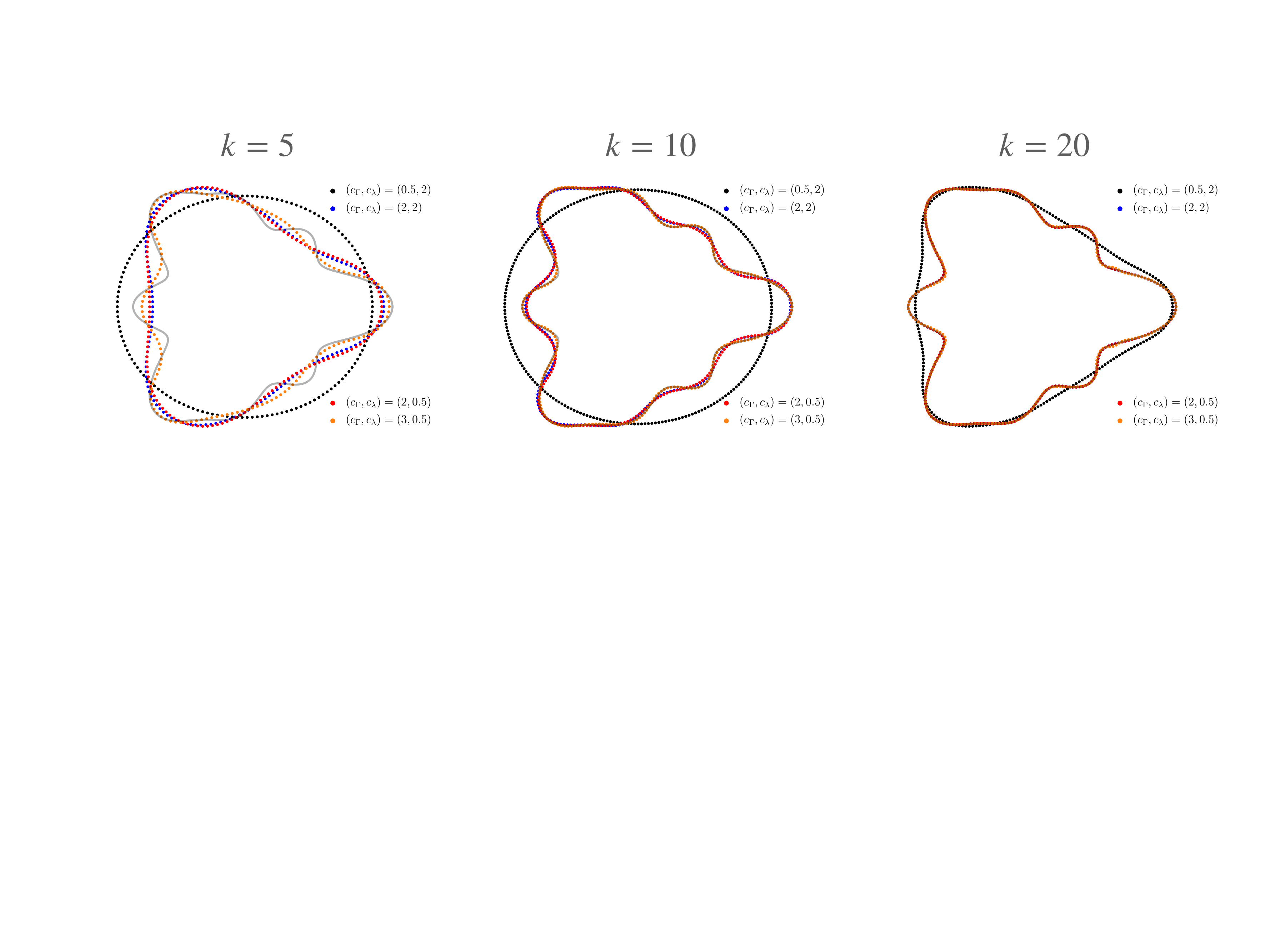}
\caption{{\bf (Section 5.1)} Comparison of the shape reconstruction at $k=5$, $10$ and $20$ using $(c_\Gamma, c_\lambda)=(0.5,2)$, $(2,2)$, $(2,0.5)$ and $(3,0.5)$ when the shape and impedance function are not known.}\label{fig:ex1c_final_shape}
\end{figure}

\begin{figure}[h]
 \center
\includegraphics[width=\linewidth]{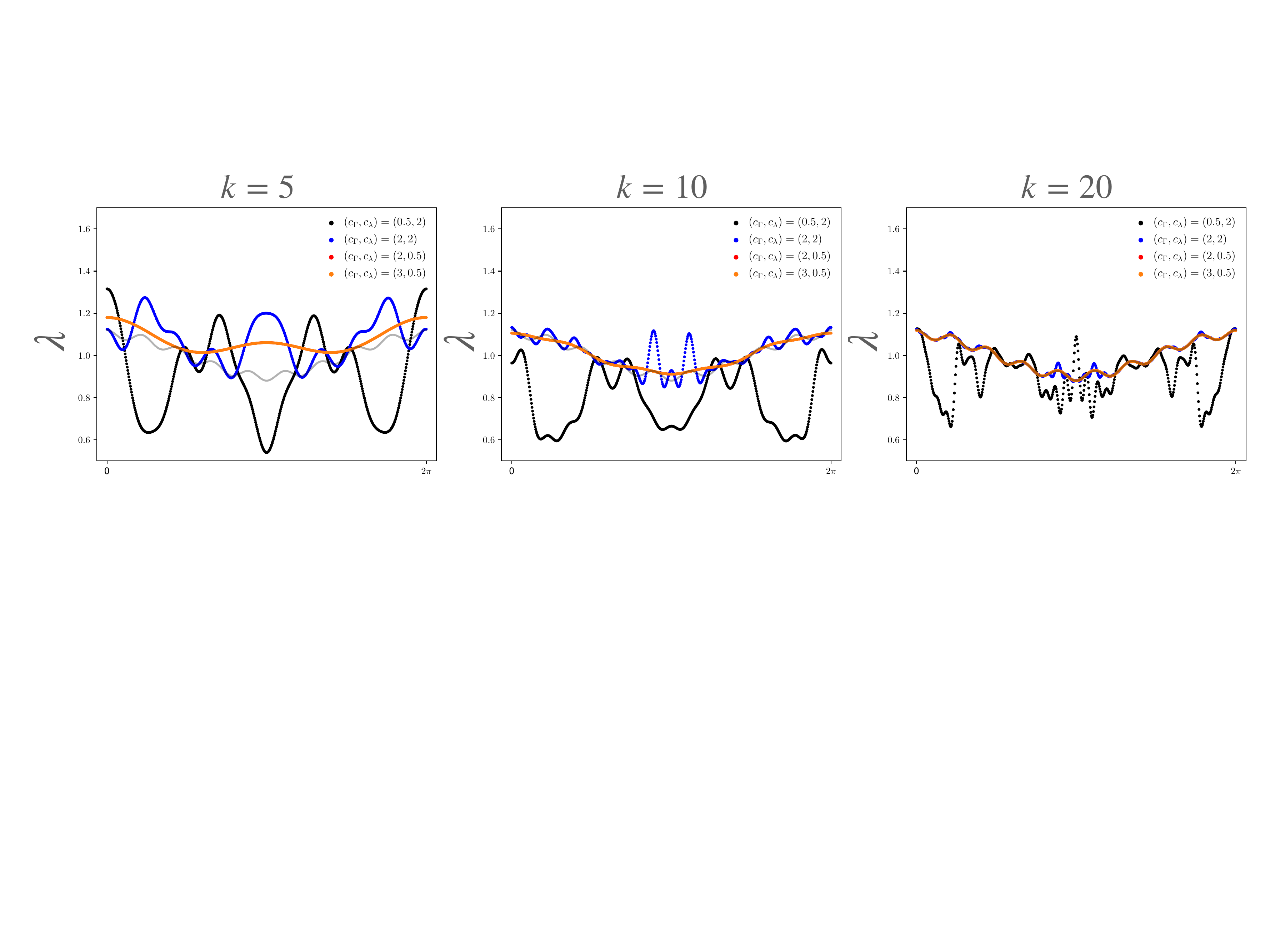}
\caption{{\bf (Section 5.1)} Comparison of the impedance function reconstruction at $k=5$, $10$ and $20$ using $(c_\Gamma, c_\lambda)=(0.5,2)$, $(2,2)$, $(2,0.5)$ and $(3,0.5)$ when the shape and impedance function are not known.}\label{fig:ex1c_final_impedance}
\end{figure}

\begin{figure}[h]
 \center
\includegraphics[width=0.7\linewidth]{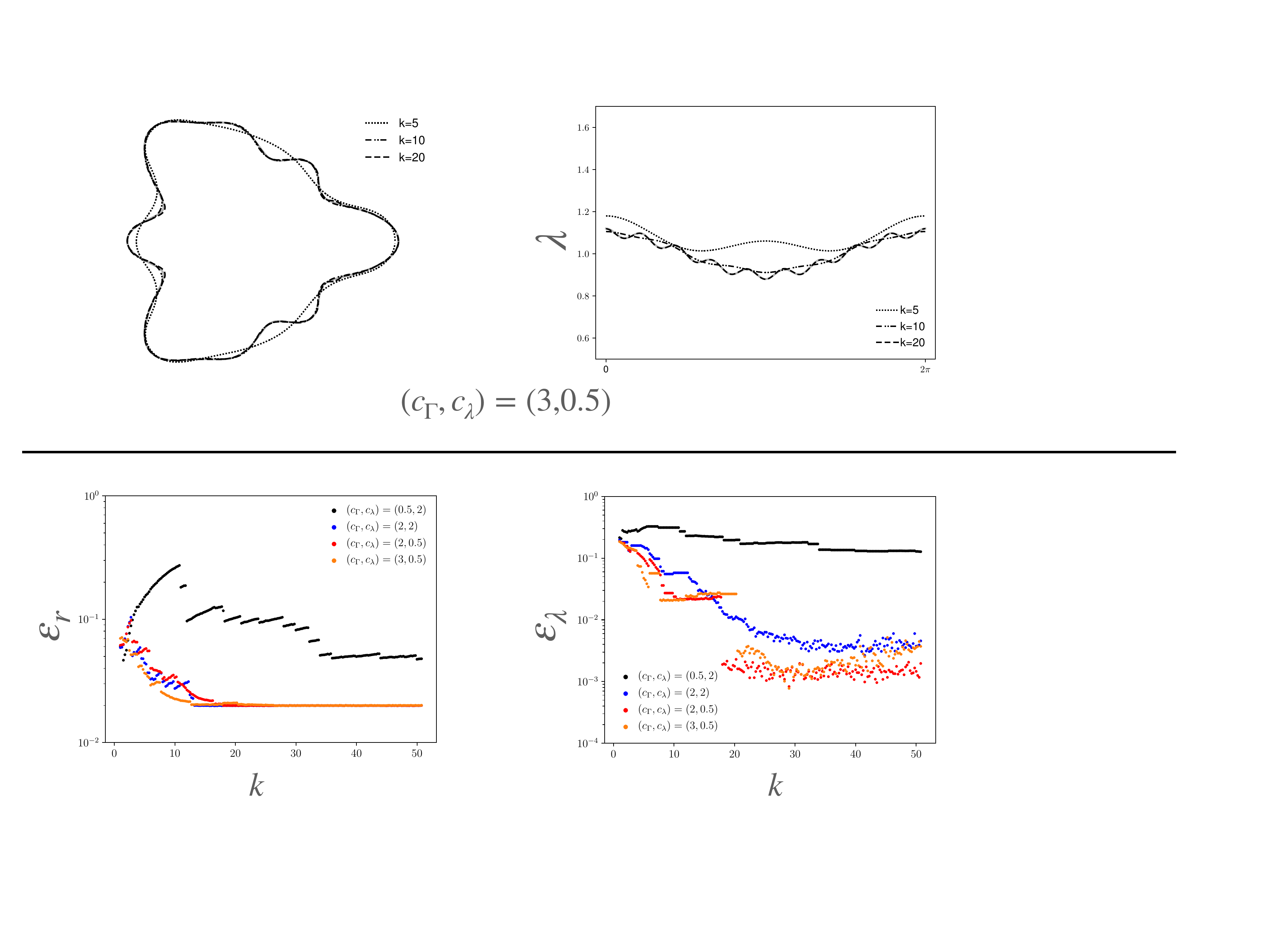}
\caption{{\bf (Section 5.1)} Relative residual and error in the reconstructed $\lambda$ using $(c_\Gamma, c_\lambda)=(0.5,2)$, $(2,2)$, $(2,0.5)$ and $(3,0.5)$ for the case where both the shape and impedance function are not known.}\label{fig:ex1c_final_residues}
\end{figure}

\begin{figure}[h]
 \center
\includegraphics[width=0.7\linewidth]{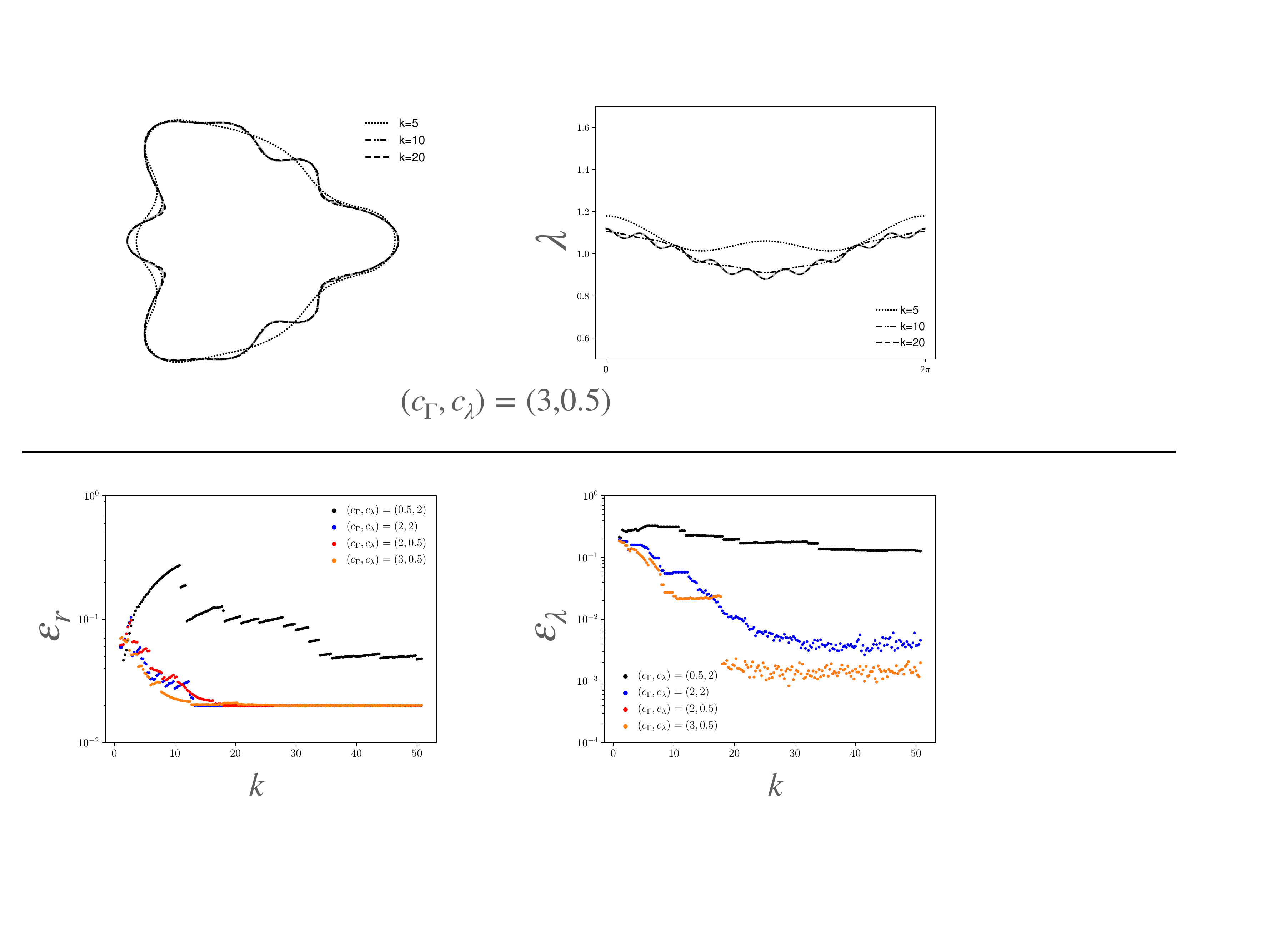}
\caption{{\bf (Section 5.1)} Reconstruction of the shape and impedance function at $k=5$, $10$ and $20$ using $(c_\Gamma, c_\lambda)=(3,0.5)$ for the case where both the shape and impedance function are not known.}\label{fig:ex1c_final_best_recons}
\end{figure}

Finally, we recover both the shape and the impedance using~\cref{alg:GN} with $(c_{\Gamma},c_{\lambda}) = (2,2)$, $(2,0.5)$, $(3,0.5)$, and $(0.5,2)$, and $\kmax=50$. In~\cref{fig:ex1c_final_shape}, we compare the reconstructions of the boundary of the obstacle for these configurations at $k=5$, $k=10$, and $k=20$. In~\cref{fig:ex1c_final_impedance}, we compare the corresponding reconstructions of the impedance. In~\cref{fig:ex1c_final_residues}, we plot the relative residual $\varepsilon_{r}$, and the error in the impedance $\varepsilon_{\lambda}$. Finally, in~\cref{fig:ex1c_final_best_recons}, we plot the reconstruction of the shape and impedance for $(c_{\Gamma},c_{\lambda})=(3,0.5)$, which was the best performing method among the four choices above. We observe high quality reconstructions for both shape and impedance for $(c_{\Gamma},c_{\lambda}) = (2,0.5), (3,0.5)$, while there are some instabilities in the reconstruction of the impedance for $(c_{\Gamma},c_{\lambda}) = (2,2)$, and stronger instabilities in the reconstruction of the impedance for $(c_{\Gamma},c_{\lambda}) = (0.5,2)$.

Note that when solving for both shape and impedance, more modes of the shape can be stably recovered as compared to the problem where the impedance is assumed to be known and fixed. The reason for this behavior can partly be attributed to the fact that fixing the impedance ends up constraining the optimization problem further. On the other hand, allowing the impedance to vary across iterates enables stable reconstruction of a larger number of modes. Moreover, the algorithm is more stable when recovering larger number of modes of the shape as compared to the impedance. A similar observation was also made in \cite{kress2001inverse}.

\subsection{Recovering the shapes of sound-soft and sound-hard scatterers \label{sec:ex2}}
\begin{figure}[h]
 \center
\includegraphics[width=\linewidth]{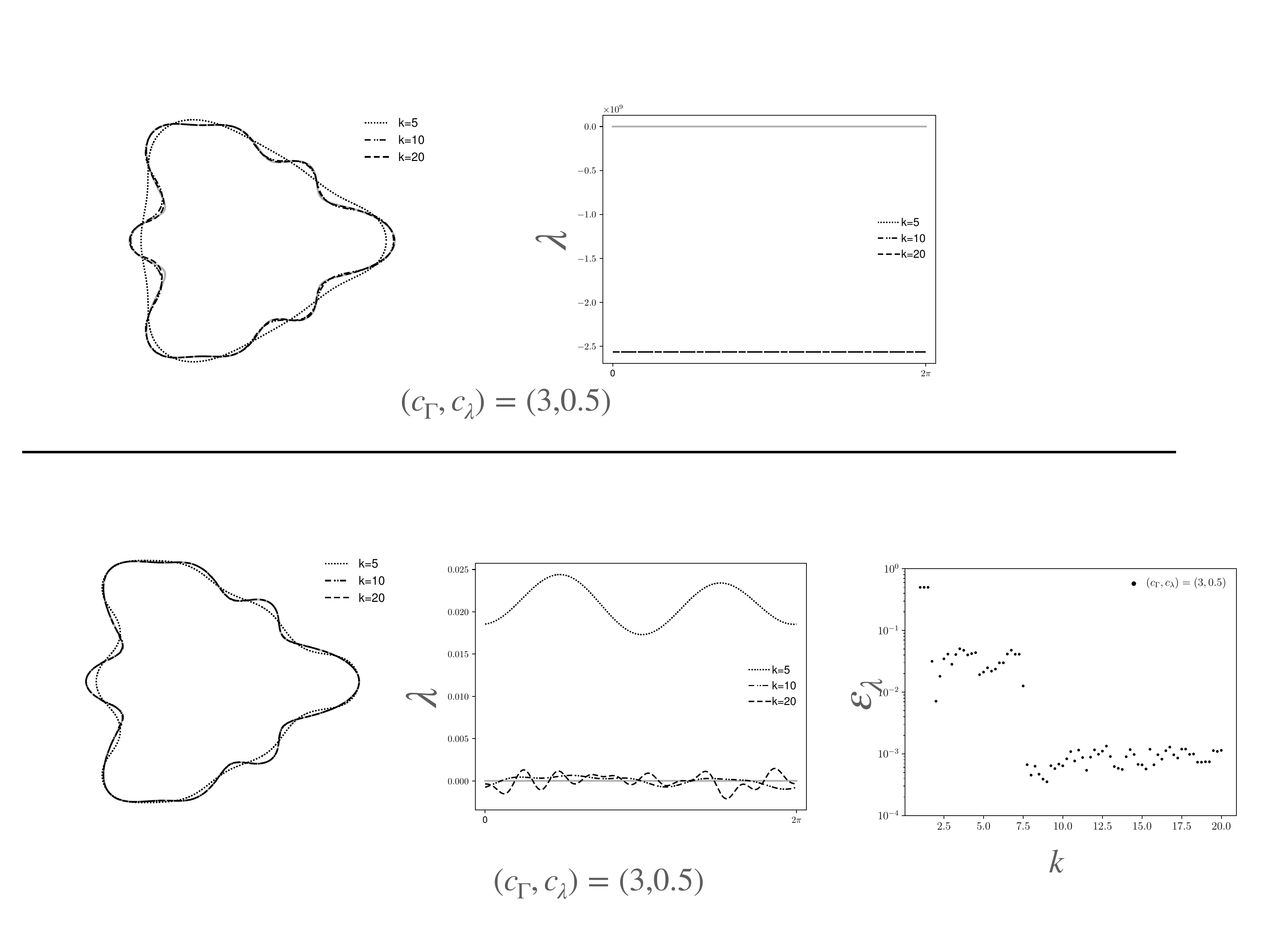}
\caption{{\bf (Section 5.2)} Shape and impedance function reconstructions at $k=5$, $10$ and $20$ with $(c_\Gamma, c_\lambda)=(3,0.5)$ when using Neumann data. We also present the error in the reconstructed $\lambda$.}\label{fig:ex2a_final}
\end{figure}

In this section, we discuss the behavior of the RLA for reconstructing shapes of sound-soft and sound-hard obstacles.   To illustrate its performance, consider the star shaped obstacle described in~\cref{sec:ex1}. For all the experiments in this section as well, we add $2\%$ noise to the scattered field measurements as discussed in~\cref{eq:noise}. First, suppose that the obstacle is a sound-hard scatterer and that the data $\ub^{\emph{meas}}$ is generated by imposing the Neumann data corresponding to the total field is $0$ on the boundary.  Sound-hard scatterers are just a special case of the scatterers with impedance boundary conditions with the impedance function set to $0$. In~\cref{fig:ex2a_final}, we plot the reconstruction of the shape of the obstacle, the impedance at $k=5,10$, and $15$, and the error in the computed impedance with $(c_{\Gamma},c_{\lambda}) = (3,0.5)$. Note that the method is able to recover an impedance that tends to zero as $k$ increases.

\begin{figure}[h]
 \center
\includegraphics[width=0.7\linewidth]{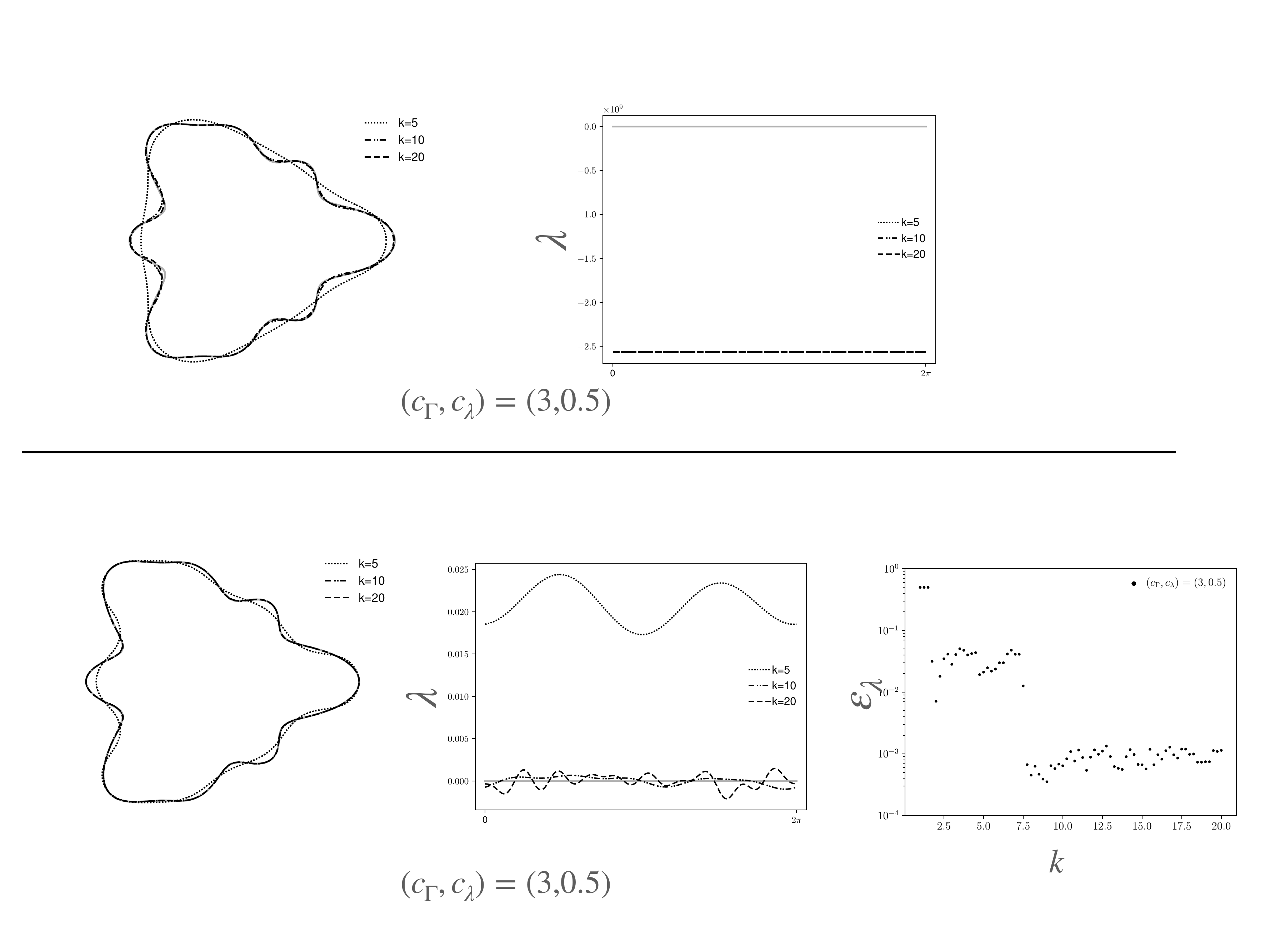}
\caption{{\bf (Section 5.2)} Shape and impedance function reconstructions at $k=5$, $10$ and $20$ with $(c_\Gamma, c_\lambda)=(3,0.5)$ when using Dirichlet data.}\label{fig:ex2b_final}
\end{figure}

Next, suppose that the obstacle is a sound-soft scatterer and that the data $\ub^{\emph{meas}}$ is generated by imposing the Dirichlet data corresponding to the total field is $0$ on the boundary. Sound-soft scatterers are also a special case of the scatterers with impedance boundary condition with the impedance function set to $\pm \infty$.
In~\cref{fig:ex2b_final}, we plot the reconstruction of the shape of the obstacle, and the impedance at $k=5,10$, and $15$, with $(c_{\Gamma},c_{\lambda}) = (3,0.5)$. 

When recovering the shape of a sound-soft scatterer, we observe the expected behavior of the impedance function diverging to $\infty$. At some Gauss-Newton iterate, the least square system comprising of the Frech\'et derivatives $\partial_{\gamma} \cF$, and $\partial_{\lambda} \cF$ is rank deficient. In fact, the Frech\'et derivative $\partial_{\lambda} \cF$ approaches the $0$ matrix. The boundary data associated with the measured scattered field for the true unknown obstacle is $u\equiv 0$ on $\Gamma$. As the reconstruction of the obstacle approaches the exact solution, 
the Frech\'et derivative $\partial_{\lambda}\cF \cdot \delta \lambda$ would be the solution to a homogeneous impedance boundary value problem independent of $\delta \lambda$ since $u\to 0$. By uniqueness of solutions to the impedance boundary value problem, we conclude that the Frech\'et derivative must be identically zero. Thus, the rank of $\partial_{\lambda} \cF$ can be used as a monitor function for determining if the object is in fact a sound-soft obstacle. Even though the impedance function diverges due to the rank deficiency of the Frech\'et derivative, the algorithm is still able to recover the shape of the obstacle to high fidelity.

\subsection{Recovering high-resolution features \label{sec:ex3}} %\label{ex:2}
In this section, we turn our attention to recovering high frequency contents of both the shape of the obstacle and the impedance function. To test this, we consider two plane like geometries which we denote by $\Omega_{1}$, and $\Omega_{2}$, where the boundary $\Omega_{2}$ has higher frequency features as compared to $\Omega_{1}$. The impedance function $\lambda:[0,2\pi]\rightarrow\mathbb{R}^2$ is given by
 \begin{equation*}
 \lambda(t)=
 \begin{cases} 
-0.1 t/\pi+0.6, \quad t\leq \pi, \\
0.1t/\pi + 0.4, \quad t>\pi.
 \end{cases}
 \end{equation*}
The impedance function is merely continuous with corners at $\pi$ and the origin, with a discontinuity in the first derivative at those points. 
The scattered field data is generated up to $\kmax = 80$. At this frequency, the domain $\Omega_{1}$ is contained in a bounding box which is $33.8\ell_{0} \times 25.3 \ell_{0}$, where $\ell_{0}=2\pi/\kmax$ is the corresponding wavelength. The perimeter of the obstacle is $258.2 \ell_{0}$. Similarly,  the domain $\Omega_{2}$ is contained in a bounding box which is $33.8\ell_{0} \times 25.1 \ell_{0}$, and it's perimeter is $365.2 \ell_{0}$.

\begin{figure}[h]
 \center
\includegraphics[width=\linewidth]{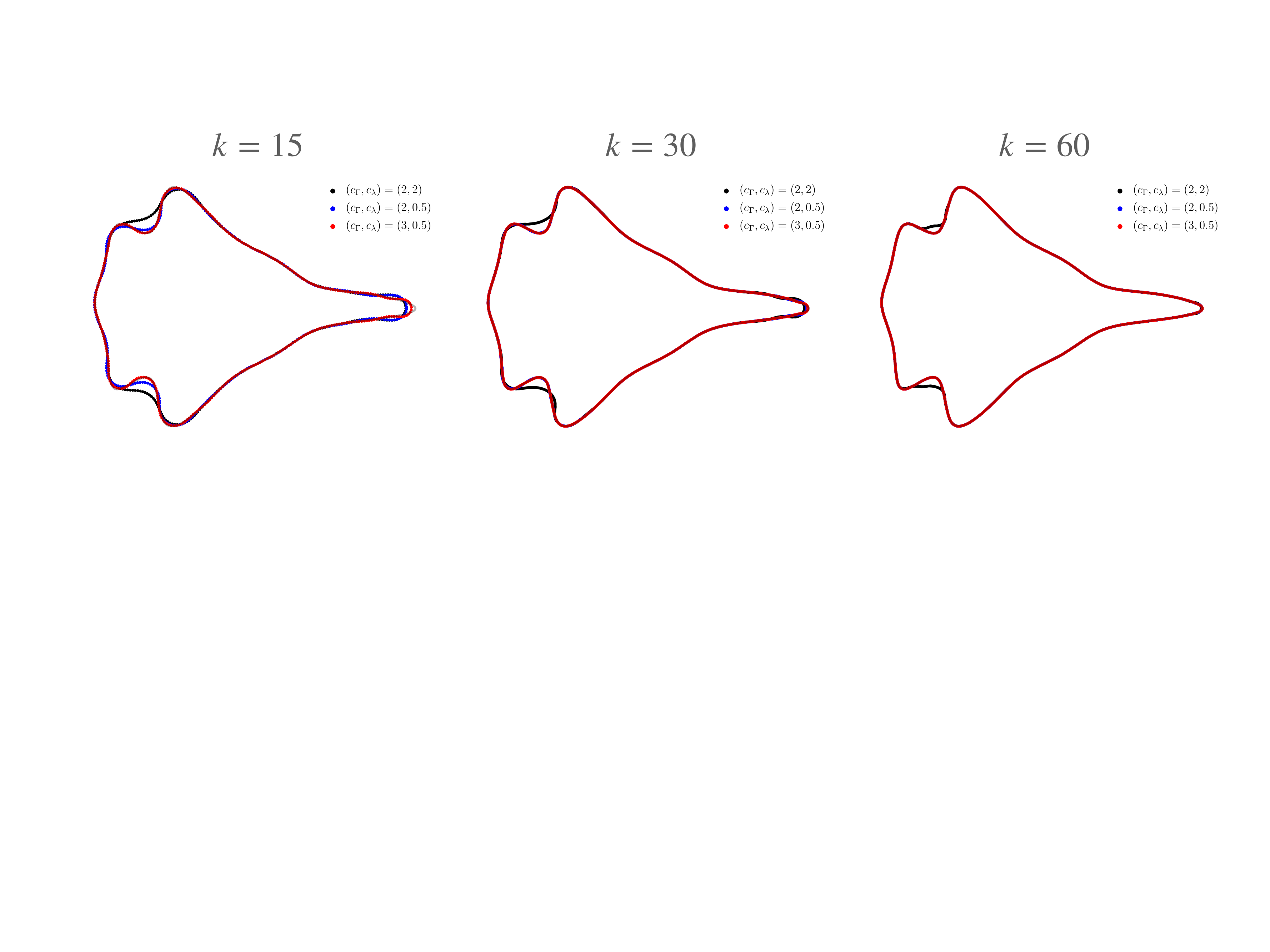}
\caption{{\bf (Section 5.3)} Comparison of the shape reconstruction at $k=15$, $30$ and $60$ for $(c_\Gamma, c_\lambda)=(2,2)$, $(2,0.5)$ and $(3,0.5)$ for the shape $\Omega_1$ and the impedance function in Section 5.3. }\label{fig:ex3a_final_shape}
\end{figure}

\begin{figure}[h]
 \center
\includegraphics[width=\linewidth]{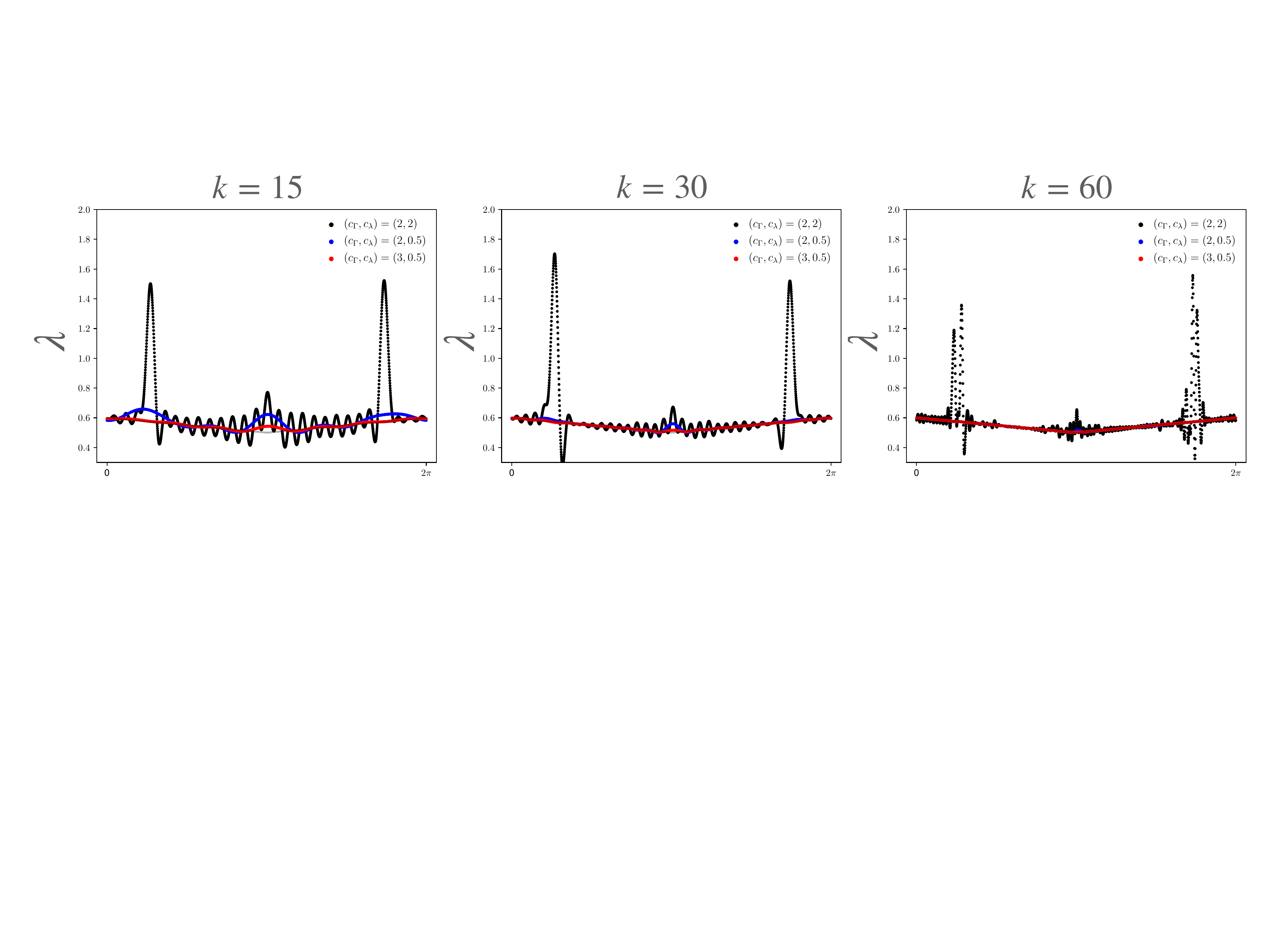}
\caption{{\bf (Section 5.3)} Comparison of the impedance function reconstruction at $k=15$, $30$ and $60$ for $(c_\Gamma, c_\lambda)=(2,2)$, $(2,0.5)$ and $(3,0.5)$ for the shape $\Omega_1$ and the impedance function in Section 5.3. }\label{fig:ex3a_final_impedance}
\end{figure}

\begin{figure}[h]
 \center
\includegraphics[width=0.7\linewidth]{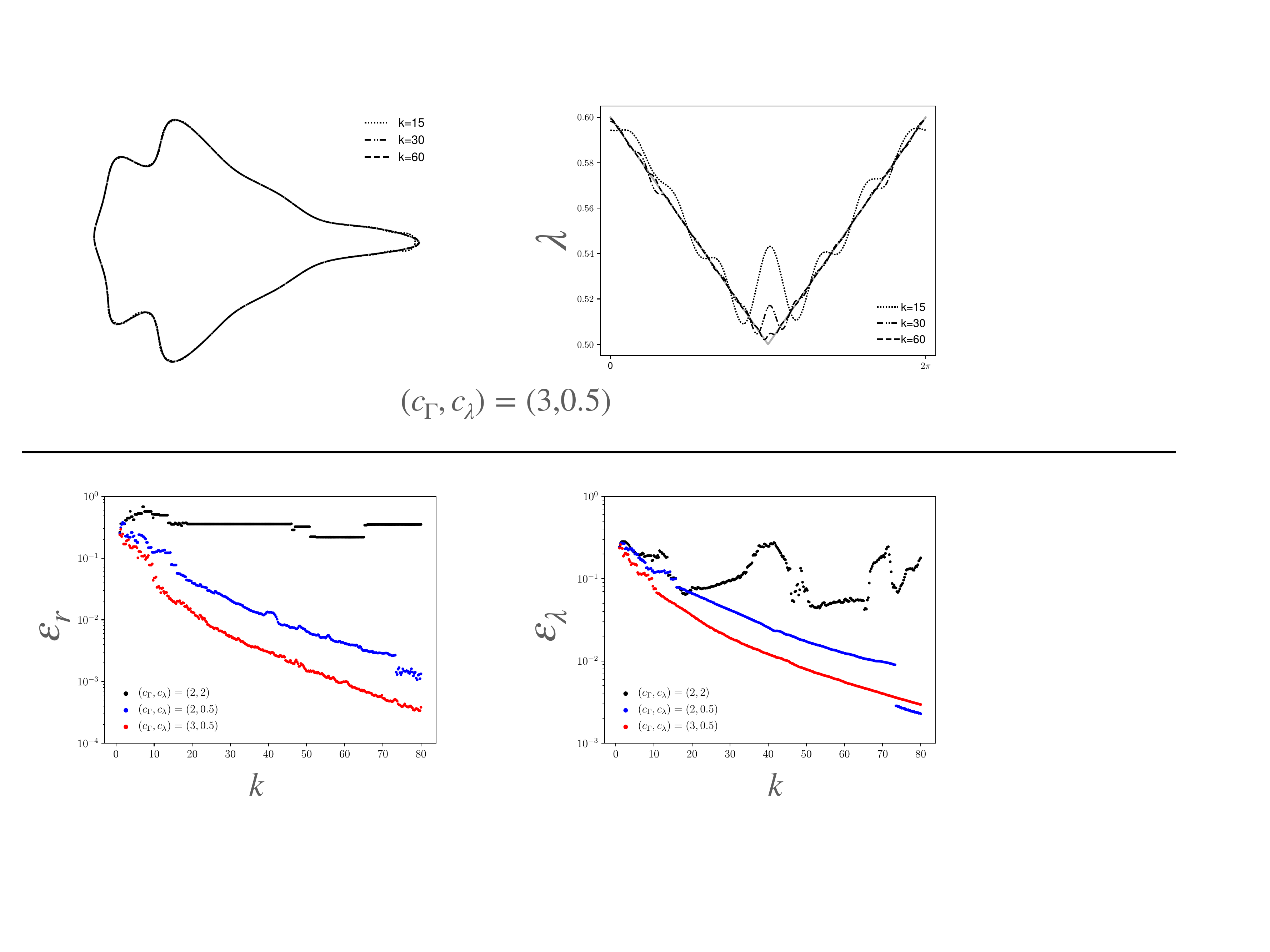}
\caption{{\bf (Section 5.3)} Relative residual and error in the reconstructed $\lambda$ for the shape $\Omega_1$ and the impedance function in Section 5.3.}\label{fig:ex3a_final_residues}
\end{figure}

\begin{figure}[h]
 \center
\includegraphics[width=0.7\linewidth]{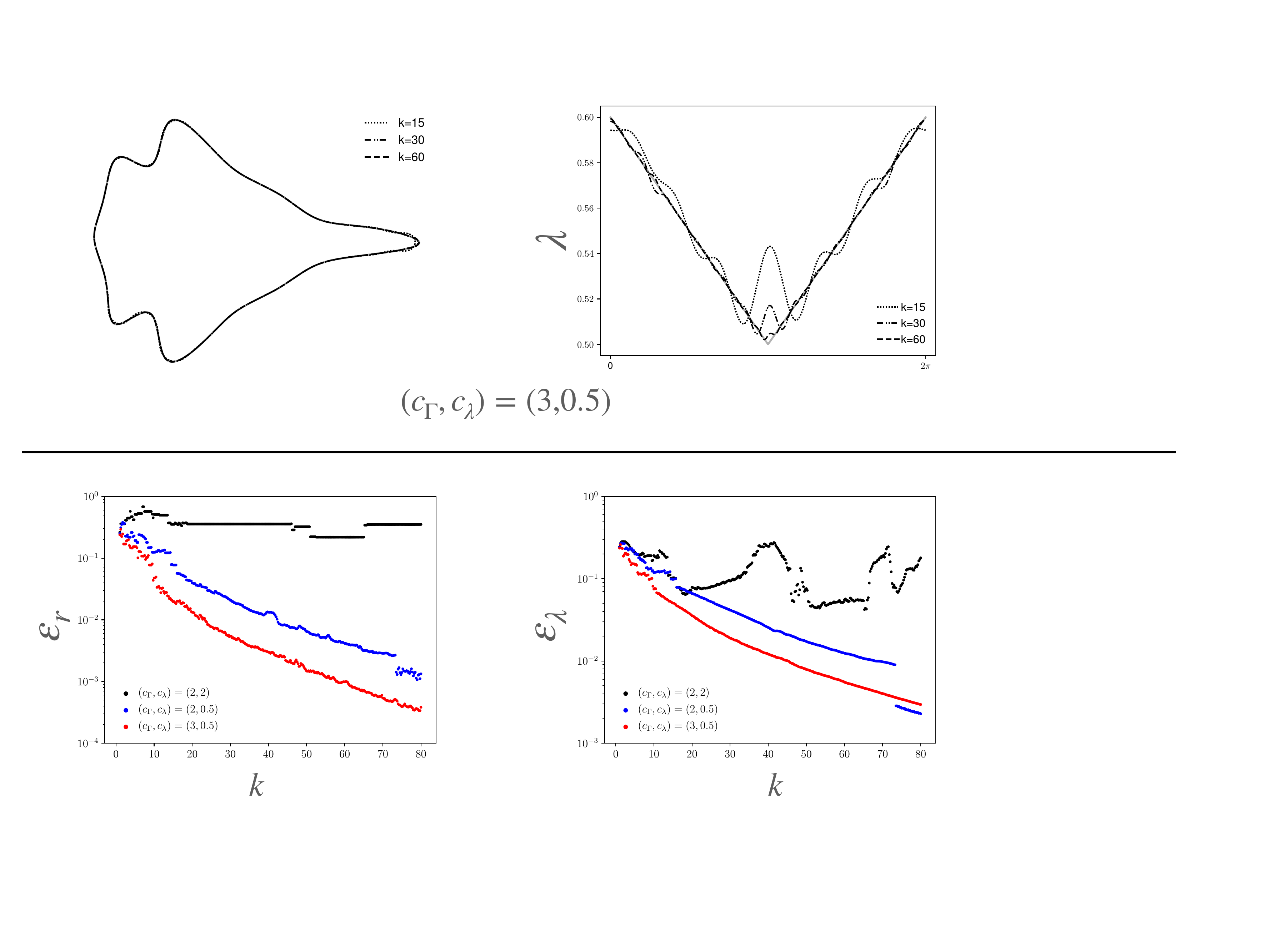}
\caption{{\bf (Section 5.3)} Best reconstruction of the shape and impedance impedance for the shape $\Omega_1$ and the impedance function in Section 5.3.}\label{fig:ex3a_final_best_recons}
\end{figure}

\begin{figure}[h]
 \center
\includegraphics[width=\linewidth]{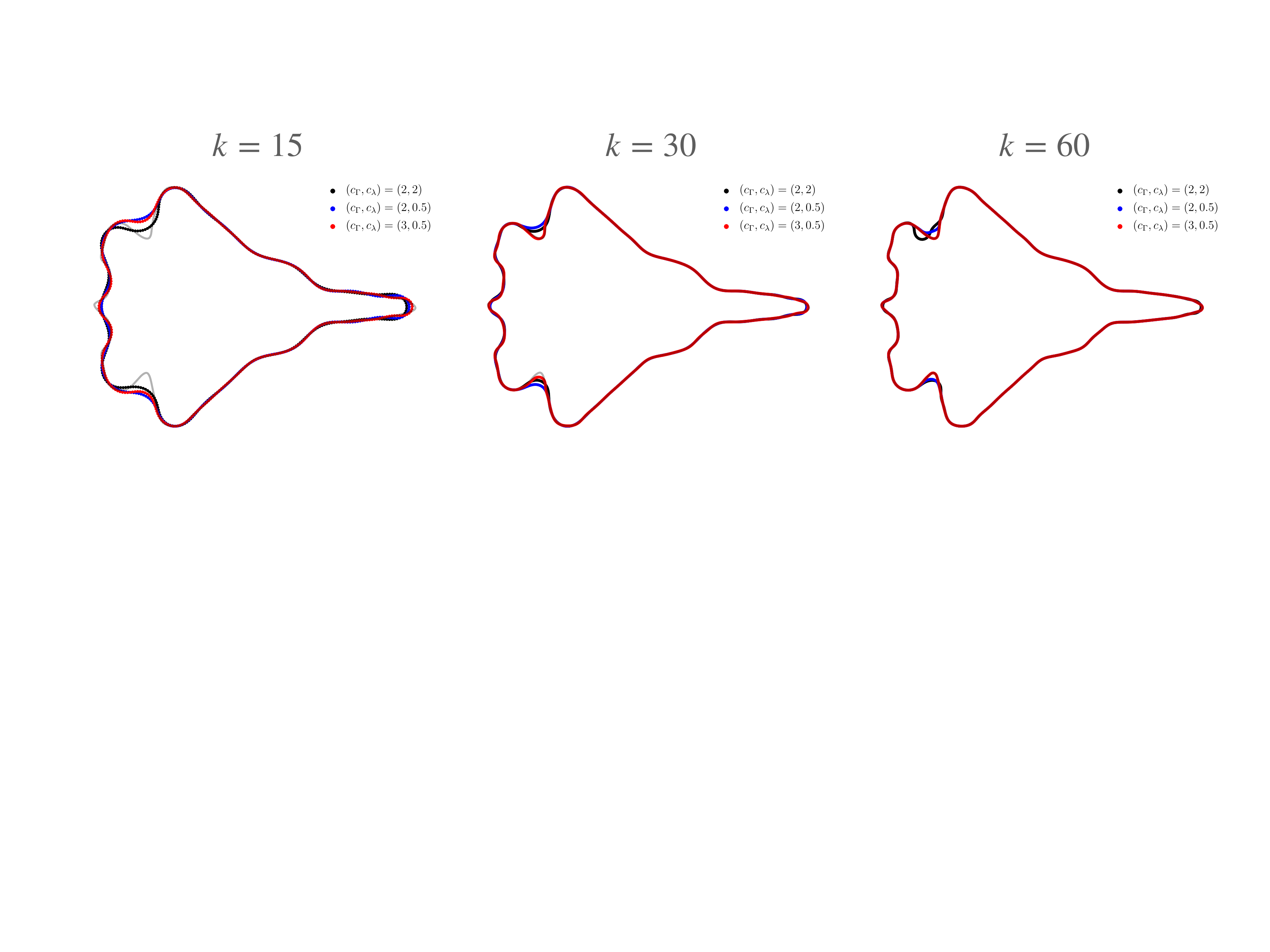}
\caption{{\bf (Section 5.3)} Comparison of the shape reconstruction at $k=15$, $30$ and $60$ for $(c_\Gamma, c_\lambda)=(2,2)$, $(2,0.5)$ and $(3,0.5)$ for the shape $\Omega_2$ and the impedance function in Section 5.3.}\label{fig:ex3b_final_shape}
\end{figure}

\begin{figure}[h]
 \center
\includegraphics[width=\linewidth]{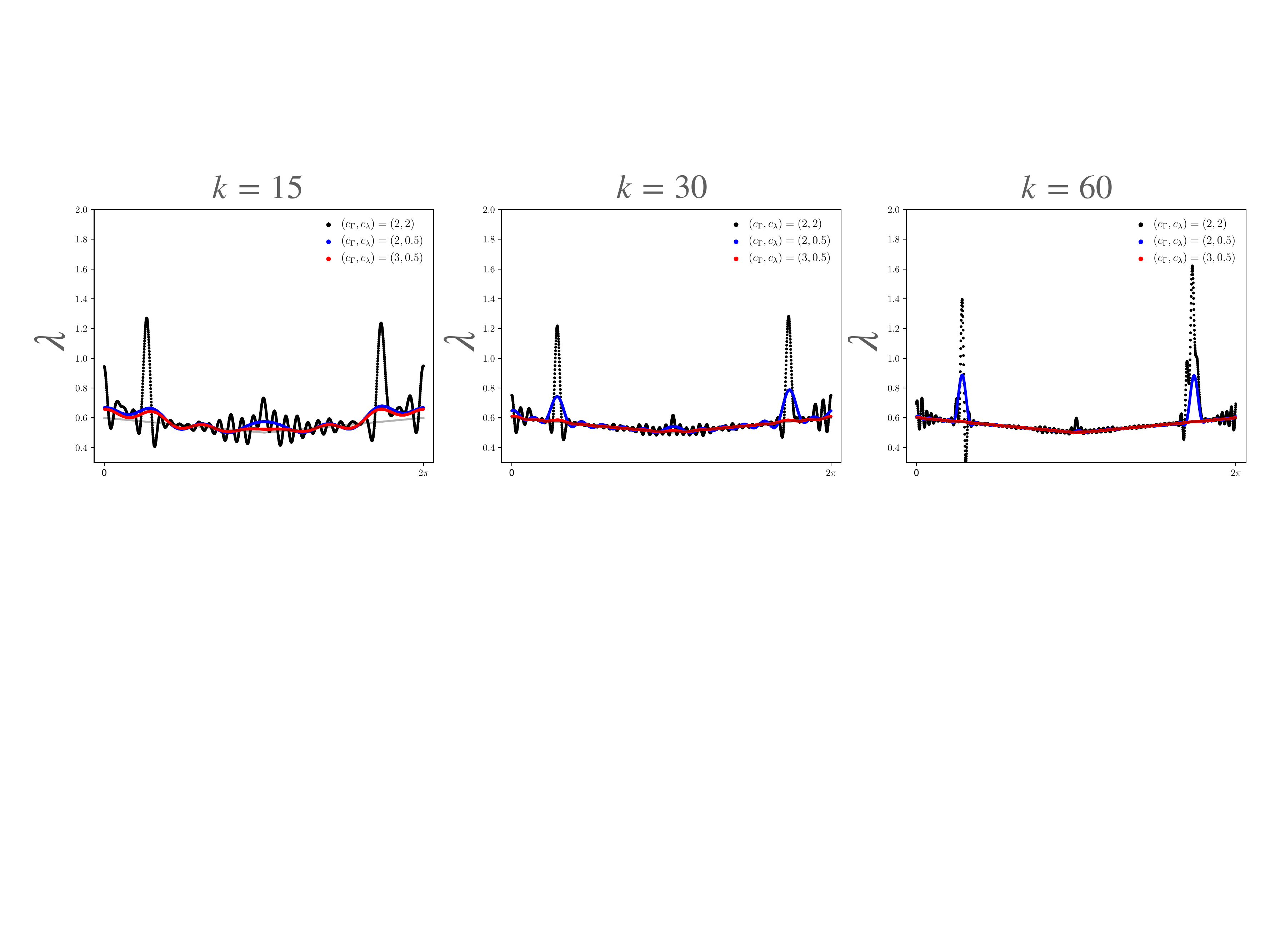}
\caption{{\bf (Section 5.3)} Comparison of the impedance function reconstruction at $k=15$, $30$ and $60$ for $(c_\Gamma, c_\lambda)=(2,2)$, $(2,0.5)$ and $(3,0.5)$ for the shape $\Omega_2$ and the impedance function in Section 5.3.}\label{fig:ex3b_final_impedance}
\end{figure}

\begin{figure}[h]
 \center
\includegraphics[width=0.7\linewidth]{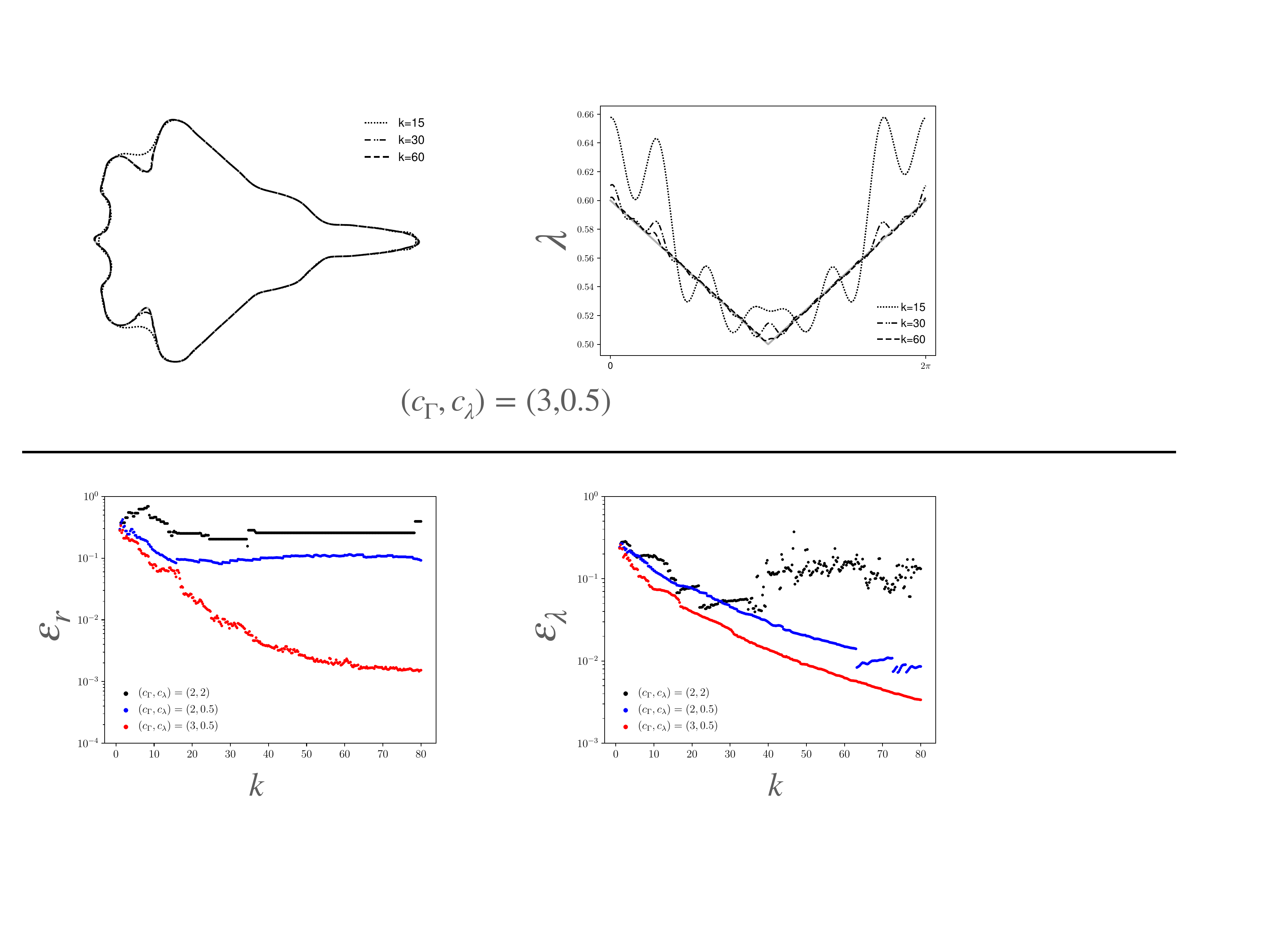}
\caption{{\bf (Section 5.3)} Relative residual and error in the reconstructed $\lambda$ for the shape $\Omega_2$ and the impedance function in Section 5.3.}\label{fig:ex3b_final_residues}
\end{figure}

\begin{figure}[h]
 \center
\includegraphics[width=0.7\linewidth]{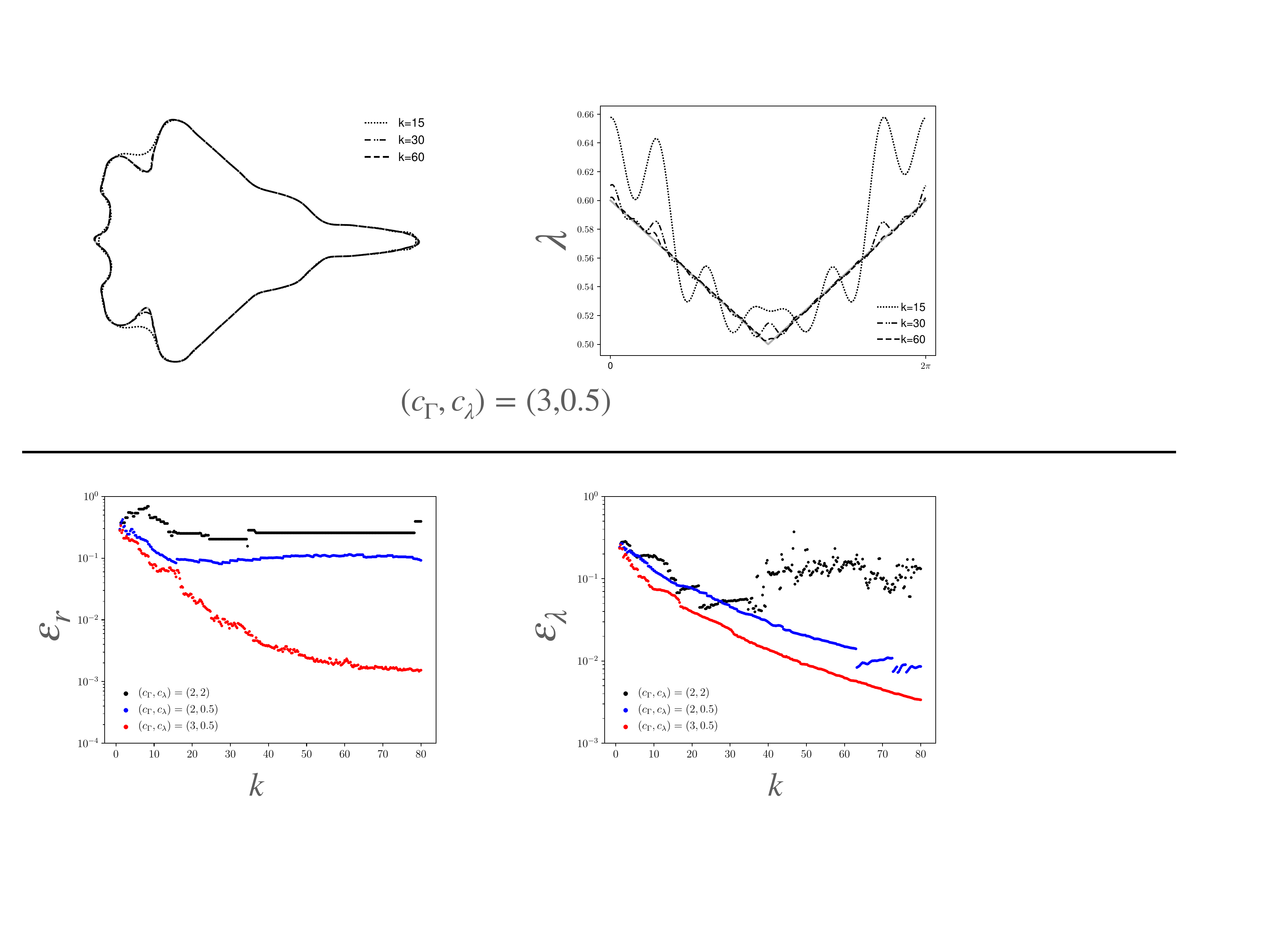}
\caption{{\bf (Section 5.3)} Best reconstruction of the shape and impedance impedance for the shape $\Omega_2$ and the impedance function in Section 5.3.}\label{fig:ex3b_final_best_recons}
\end{figure}

In~\cref{fig:ex3a_final_shape}, we compare the reconstructions of the boundary of the obstacle $\Omega_{1}$ at $k=15$, $k=30$, and $k=60$, with $(c_{\Gamma},c_{\lambda})=(2,2),(2,0.5)$, and $(3,0.5)$. In~\cref{fig:ex3a_final_impedance}, we compare the corresponding reconstructions of the impedance. In~\cref{fig:ex3a_final_residues}, we plot the relative residual $\varepsilon_{r}$, and the error in impedance $\varepsilon_{\lambda}$, and finally in~\cref{fig:ex3a_final_best_recons}, we plot the reconstruction of the shape and impedance for $(c_{\Gamma},c_{\lambda})=(3,0.5)$, which was the best performing method among the three choices above.~\Cref{fig:ex3b_final_shape,fig:ex3b_final_impedance,fig:ex3b_final_residues,fig:ex3b_final_best_recons} are the corresponding results for $\Omega_{2}$.

We observe high quality reconstructions for both shape and impedance for $(c_{\Gamma},c_{\lambda}) = (3,0.5)$ for both $\Omega_{1}$, and $\Omega_{2}$. With $(c_{\Gamma},c_{\lambda}) = (2,0.5)$, we observe high quality reconstructions for both the shape and impedance for $\Omega_{1}$, however the algorithm in unable to recover the deeper cavities in $\Omega_{2}$. This results in instabilities in the reconstruction of the impedance in the vicinity of the cavities. As before for the example in~\cref{sec:ex1}, we are unable to stably recover the shape or impedance for either $\Omega_{1}$ or $\Omega_{2}$ with $(c_{\Gamma},c_{\lambda})=(2,2)$.

\section{Open problems} \label{s:openproblems}
While the RLA is very effective for the solution of inverse obstacle scattering problems, the behavior of the residue $\varepsilon_{r}$ for various configurations point to several open questions that require further investigation. The lack of improvement in reconstruction at some of the frequencies, particularly for $k\leq 10$ in~\cref{fig:ex1a_final} seems to indicate that the Gauss-Newton algorithm for the single frequency minimization problem is not converging to a local minimum for those frequencies. We can show that the result does not converge to a global minimum of the single frequency problem by using a different optimization method for minimizing the single frequency problem. A convenient choice for this purpose is to use a damped Gauss-Newton algorithm where the step size is chosen to be $\min(c_{s}\pi/k, 1)$, with $c_{s} = 1,2$. In a slight abuse of notation, let $c_{s} = 0$ denote the case where no step size control is implemented. The intuition for the particular scaling of the step size control stems from the fact that we expect the objective function $\| \ub^{\emph{meas}} - \cF_{k} \|$ to oscillate at the wavelength of the underlying Helmholtz problem which is $2\pi /k$. 

\begin{figure}[h]
 \center
\includegraphics[width=\linewidth]{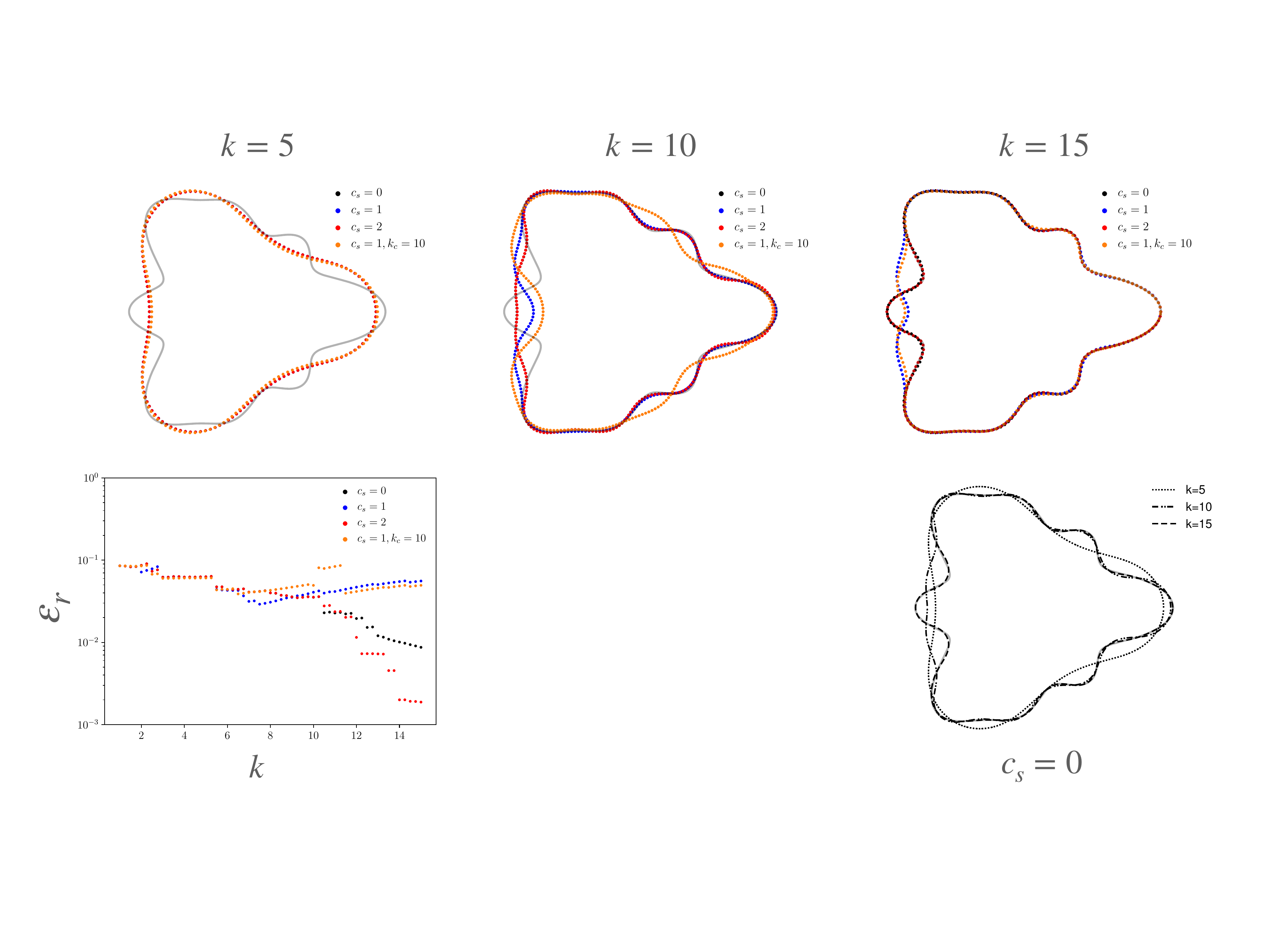}
\caption{{\bf (Section 6)} Shape reconstructions using the control steps $c_{s}=0,1$ and $2$ at $k=5,10$, $15$. We also present the shape reconstruction using multiple frequency data, $k_c=10$, with step size control $c_s=1$. The relative residues $\varepsilon_{r}$, and the reconstruction as a function of frequency for no step size control. In these experiments, we assume the impedance function is known.}\label{fig:ex_openprob12_final}
\end{figure}

Consider the recovery of the shape of the obstacle assuming the impedance known where the scattered field measurements are generated using the domain $\Omega$ and the impedance $\lambda$ defined in~\cref{sec:ex1}. We use $c_{\Gamma}=2$ for all the reconstructions. In~\cref{fig:ex_openprob12_final}, we plot the reconstructions using $c_{s}=0,1$ and $2$ at $k=5,10$, $15$, along with the relative residues $\varepsilon_{r}$, and the reconstruction as a function of frequency for no step size control. 

Referring to~\cref{fig:ex_openprob12_final}, note that for $6\leq k\leq 8$, $\varepsilon_{r}$ for $c_{s}=1$ is lower as compared to $c_{s}=0$ or $c_{s}=2$. This indicates that~\cref{alg:GN} (corresponding to $c_s=0$) has not converged to a global minimum of the loss function. However, the final reconstruction obtained using $c_{s}=1$ is much poorer as compared to $c_{s}=0$, or $c_{s}=2$, due to the solution being stuck in a local minimum at low frequencies. The termination of the single frequency optimization problem due to an increase in residue seems to play an important role in indicating if the reconstruction at that frequency is moving out of the global basin of attraction of the exact solution. 
Thus, one might need to strike a balance between finding the global minimum of the single frequency optimization problem with trying to stay in the basin of attraction for the multifrequency problem.

We investigate this issue further and examine if this issue is a result of the reduction of the multifrequency problem to a collection of single frequency problems. To answer this question, we use the cumulative loss function given by~\cref{eq:single_freq_inv_cum}, for $k\leq k_{c}$ and then switch to a sequence of single frequency minimization problems. While this improves the final reconstruction marginally as compared to using the single frequency loss function, the solution still seems to be stuck in a similar local minimum. Thus, the above behavior is not due to the reduction of the multifrequency problem to a sequence of single frequency problems. An outstanding question is: what optimization algorithm and what termination criterion/monitor functions should be used for minimizing  $\| \ub^{\emph{meas}} - \cF_{k} \|$, such that the sequence of minimizers in $k$ stay in the basin of attraction of the true minimizer for the multifrequency problem?

It is unclear from these experiments whether the issue lies in the choice of the optimization method or the optimization problem itself. There are several open problems along this line of investigation, particularly for non star-shaped domains. Under what hypotheses, does the single frequency optimization problem have a unique minimizer which can be stably recovered in finite precision arithmetic? In particular, what constraints need to be put on the spaces/subsets for recovering the shape of the obstacle and the impedance? How many measurements of the scattered field are required? What is the impact of noise in measurements of the scattered data on the reconstruction? Similar questions need to be addressed for the multifrequency problem as well. An additional consideration would be whether a sequence of minimizers for the single frequency problem converge in some appropriate sense to the obstacle and impedance from which the data is generated. The last question is particularly essential for an approach like the RLA to work effectively for such inverse problems.
The ability to separate out the well-posedness of the optimization problem and its behavior in the multifrequency regime will help tremendously in the construction of robust algorithms for inverse obstacle problems.

\section{Conclusions}\label{s:conclusions}
In this paper, we presented an extension of the RLA for the solution of the inverse scattering problem of recovering the shape and impedance boundary function of an impenetrable obstacle using multifrequency measurements of the scattered field. In this approach,  the multifrequency inverse problem is reduced to the solution of a sequence of constrained single frequency inverse problems with increasing wavenumber, wherein each single frequency inverse problem is optimized using a Gauss-Newton method with a bandlimited representation of the variables.

Using this approach, we were able to obtain high fidelity reconstruction of obstacles for which both the shape and the impedance functions had high frequency features. Similar to the conclusions made in~\cite{kress2001inverse}, we observe that when recovering the shape and impedance function of an obstacle simultaneously, it is more beneficial to recover more modes of the shape than the impedance function. Moreover, using this approach, we were also able to recover the shape of sound-soft, and sound-hard obstacles. When recovering sound-soft obstacles, the rank of the Frech\'et derivative with respect to the impedance can be used as a monitor function for detecting whether the scattered data was generated by a sound-soft obstacle.

Even though the RLA performed well in recovering the shape and impedance in the examples presented, there are several open questions regarding this approach. Is it possible to guarantee that the algorithm converges to a global minimum? What is the best way to constrain the shape of the obstacle? What is the effect of the noise to the multifrequency scheme?

The approach outlined in this paper extends almost immediately to the case of higher-order impedance operators. Impedance boundary conditions are also often used for approximating thin surface coatings where the impedance function is proportional to the depth of the coating. A natural extension would be to recover the depth of the coating from far-field measurements of the scattered field. Impedance boundary conditions/generalized impedance boundary conditions are also often used in applications where the unknown obstacle has some dissipation and a different but constant wave speed from its surrounding medium. In this setup, one could compare the solution of the inverse problem using one of three approaches: recovering the unknown sound speed (which in this case is discontinuous); solve an inverse obstacle problem to recover the shape and impedance or generalized impedance; and finally solve an inverse obstacle problem where the sound speed inside the obstacle is treated as an additional unknown. This line of inquiry along with the open problems presented in~\cref{s:openproblems} are currently being vigorously pursued and will be reported at a later date.

Finally, while most of the RLA extends naturally to the recovery of obstacles in three dimensions, handling the geometry of space of surfaces is particularly challenging. The ability to specify and update complicated bandlimited surfaces in the optimization loop will play a critical role for many inverse obstacle scattering problems in three dimensions.

\section{Acknowledgments}
The authors would like to thank Leslie Greengard for many useful discussions.

\FloatBarrier

\bibliography{./Biblio}

\end{document}